\documentclass[9pt,twocolumn,twoside]{IEEEtran}

\renewcommand{\tilde}{\widetilde}
\renewcommand{\hat}{\widehat}
\newcommand{\cd}{(\cdot)}

\usepackage{color}
\usepackage{dsfont}
\usepackage{graphics} 
\usepackage{epsfig} 
\usepackage{amsmath} 
\usepackage{amssymb}  
\usepackage{subfigure}
\usepackage{stfloats}
\usepackage{array}
\newcommand{\one}{\mathds{1}}
\newcommand{\ad}{&\!\!\!\disp}
\newcommand{\aad}{&\disp}
\newcommand{\bed}{\begin{displaymath}}
\newcommand{\eed}{\end{displaymath}}
\newtheorem{Theorem}{Theorem}[section]
\newtheorem{Lemma}{Lemma}[section]
\newtheorem{Definition}{Definition}[section]
\newtheorem{Remark}{Remark}[section]
\newtheorem{Corollary}{Corollary}[section]
\newcounter{MYtempeqncnt}

\begin{document}


\title{Adaptive Search Algorithms for Discrete Stochastic Optimization: A Smooth Best-Response Approach}

\author{Omid~Namvar~Gharehshiran,
        Vikram~Krishnamurthy,
        and George~Yin.
\thanks{O. N. Gharehshiran and V. Krishnamurthy are with the department of Electrical and Computer
Engineering, University of British Columbia, Vancouver, V6T 1Z4, Canada (e-mail: omidn@ece.ubc.ca; vikramk@ece.ubc.ca).
This research was supported by the NSERC Strategic grant and the Canada Research Chairs program.}
\thanks{G. Yin is with the Department of Mathematics, Wayne State University, Detroit, MI 48202, USA (e-mail: gyin@math.wayne.edu).
This research
 was supported in part by the Army Research Office under grant
W911NF-12-1-0223.}%
}


\maketitle
\begin{abstract}
This paper considers simulation-based optimization of the performance of a regime-switching stochastic system over a finite set of feasible configurations. Inspired by the stochastic fictitious play learning rules in game theory, we propose an adaptive simulation-based search algorithm that uses a smooth best-response sampling strategy and tracks the set of global optima, yet distributes the search so that most of the effort is spent on simulating the system performance at the global optima. The algorithm converges weakly to the set of global optima even when the observation data is correlated (as long as a weak law of large numbers holds). Numerical examples show that the proposed scheme yields a faster convergence for finite sample lengths compared with several existing random search and pure exploration methods in the literature.
\end{abstract}


\def\x{\boldsymbol{x}}
\def\L{\mathcal{L}}
\def\F{\mathcal{F}}
\def\C{\mathcal{N}}
\def\SS{\mathcal{S}}
\def\X{\boldsymbol{X}}
\def\T{\mathcal{T}}
\def\H{\mathcal{H}}

\def\A{\mathcal{A}}
\def\K{\mathcal{K}}

\def\U{\mathcal{U}}
\def\a{\boldsymbol{a}}
\def\t{\mathbf{t}}
\def\e{\mathbf{e}}
\def\z{\boldsymbol{z}}
\def\p{\mathbf{p}}
\def\h{\boldsymbol{h}}
\def\R{\boldsymbol{R}}
\def\RR{{\mathbb{R}}}
\def\RRk{\RR^r}
\newcommand{\disp}{\displaystyle}
\def\p{\boldsymbol{p}}
\newcommand{\br}{\boldsymbol{b}^{\gamma}}
\newcommand{\bbr}{b^{\gamma}}
\def\f{\tilde{\boldsymbol{f}}}

\newcommand{\beq}[1]{\begin{equation} \label{#1}}
\newcommand{\eeq}{\end{equation}}
\newcommand{\barray}{\begin{array}{c}}
\newcommand{\earray}{\end{array}}
\renewcommand{\bar}{\overline}
\newcommand{\argmin}{\operatornamewithlimits{argmin}}
\newcommand{\argmax}{\operatornamewithlimits{argmax}}

\newcommand{\bea}{$$\begin{array}{rl}}
\newcommand{\eea}{\end{array}$$}

\def\rk{r_{ij}^k}
\def\drk{\bar{r}_{ij}^k}
\def\pk{p^k}
\def\Xk{X^k}
\def\Nkc{\mathcal{N}_k^c}
\def\Nk{\mathcal{N}_k}
\def\M{\mathcal{M}}
\def\vare{\varepsilon}
\def\rr{\mathbb{R}}
\def\go{\mathcal{S}}
\newcommand{\ee}{{\mathbb E}}
\def\snk{s_n^k}
\def\s{\boldsymbol{s}}

\def\lb{\left[}
\def\rb{\right]}
\def\lbr{\left\lbrace}
\def\rbr{\right\rbrace}

\section{Introduction}
\label{sec:introduction}

\IEEEPARstart{D}{iscrete} stochastic optimization problems arise in operations research~\cite{norkin1998optimal,swisher2000survey},
manufacturing engineering~\cite{park2005line}, and
communication networks~\cite{krishnamurthy2004spreading,berenguer2005adaptive}.
These problems are intrinsically more difficult to solve than their deterministic counterparts due to the non-availability of an explicit relation between the objective function and the underlying decision variables. It is therefore necessary to use stochastic simulation to estimate the objective function in such problems.

\subsection{The Problem}
\label{sec:intro-problem}
The simplest setting of a discrete stochastic optimization problem is as follows: Estimate
\beq{eq:SOP_stat}
\mathcal{S} := \argmin_{s\in\M}
F(s) = \argmin_{s\in\M} \ee \lbr f_n(s)\rbr,
\eeq
where the search space $\M = \lbr1,2,\ldots,S\rbr$ is finite, $\lbr f_n(s)\rbr$ for each $s\in\M$ is a sequence of i.i.d. random variables with finite variance but
unknown distribution, $\ee$ denotes expectation with respect to the distribution of $f_n(s)$, and $F:\M\to \RR$ is deterministic. Typically, $F\cd$ represents the expected performance of a stochastic system. Since the distribution of $\lbr f_n(s)\rbr$ is unknown, $F(s)$ cannot be evaluated analytically.

A brute force method of solving~(\ref{eq:SOP_stat}) involves an exhaustive enumeration: For each $s\in\M$ compute
$\hat{F}_N(s) = \frac{1}{N} \sum_{n=1}^N f_n(s)$ via simulation for large $N$. Then, pick $\hat{s}^* = \argmin_{s \in \mathcal{M}} \hat{F}_N(s)$. Since $\lbr f_n(s)\rbr$ for each $s \in \M$ is an i.i.d. sequence of random variables, Kolmogorov's strong law of large numbers implies that $\hat{F}_N(s) \to \ee \lbr f_n(s)\rbr$ almost surely as $N \to \infty$. This, together with the finiteness of $\M$ implies that as $N \to \infty$,
\begin{equation}
\argmin_{s\in\M} \hat{F}_N(s) \to \argmin_{s\in\M} \ee\lbr f_n(s)\rbr \;\; \textmd{w.p.1}.
\nonumber
\end{equation}
This requires $f_n(s)$ to be evaluated for each $s\in\M$ at each sampling period, and is highly inefficient since the evaluations $f_n(s)$, for $s\notin \go$, do not contribute to finding $\go$ and are wasted. The main idea here is to develop a novel adaptive search scheme that is both \emph{attracted} to the global optima $\go$ and \emph{efficient}, in the sense that it spends most of its effort simulating $\go$~\cite[Chapter 5.3]{pflug1996optimization}.

Problem~(\ref{eq:SOP_stat}) is static in the sense that the set of global minima $\go$ does not evolve with time. In this paper, we consider two extensions of the above problem: First, we allow for $\lbr f_n(s)\rbr$ to be a correlated sequence, as long as it satisfies the weak law of large numbers. Second, we
solve an adaptive variant of this problem where the set of global optima evolves with time according to the sample path of a finite-state Markov chain $\lbr\theta(n)\rbr$ with state space $\mathcal{Q} = \lbr 1,2,\ldots,\Theta\rbr$. More precisely, consider a simulation-based discrete stochastic optimization problem of the form
\beq{eq:SOP}
\go(\theta(n)) := \argmin_{s\in\M}
F(s,\theta(n)) =\argmin_{s\in\M} \ee \lbr f_n(s,\theta(n))\rbr.
\eeq
%
We assume that the Markov chain $\lbr \theta(n)\rbr$ cannot be observed and its dynamics are unknown. However, for any choice of $s \in\M$, the samples $f_n(s,\theta(n))$ can be generated via simulation.
 We further allow time correlation of the simulation data $f_n(s,\bar{\theta})$ for each $\bar{\theta}\in\mathcal{Q}$, that is more realistic in practice.

The Markov chain $\lbr \theta(n)\rbr$ in~(\ref{eq:SOP}) constitutes the so-called {\em hypermodel}~\cite{benveniste1990adaptive} for the underlying dynamics. It
represents the jump changes in the profile of the stochastic events in the system or the objective function or both. Such problems arise in a broad range of  practical applications where the goal is to track the optimal operating configuration of a stochastic system subject to time inhomogeneity. We assume that the transition probability matrix of the Markov chain $\lbr\theta(n)\rbr$ is ``close'' to the identity matrix. That is, the Markov chain has transition matrix $I+\vare Q$, where $\vare$ is a small parameter. We will refer to such a Markov chain with infrequent jumps as {\em slow Markov chain}, for simplicity. The global optima $\go(\theta(n))$ thus varies with time according to the slow Markov chain.
 In what follows,
 we refer to the above problem as ``regime-switching discrete stochastic optimization''. Tracking such time-varying sets lies at the very heart of applications of adaptive stochastic approximation algorithms.

\paragraph*{Example} Consider the problem of optimizing buffer sizes in a queueing network comprising multiple stations with buffers. Such a network may represent an assembly line in the manufacturing industry, networked-processors in parallel computing, or a communication network. Let $s$ and $\lbr X_n(s,\theta(n))\rbr$ denote the vector of buffer sizes and the sequence of random vector of service times at different stations, respectively.
The distribution of service times may jump change due to the changes in the nature of the offered services. The performance of such a system $f(s, X_n(s,\theta(n)))$ is a function of both $s$ and $\lbr X_n(s,\theta(n))\rbr$ and is often evaluated by the amortized cost of buffers minus the revenues due to the processing speed. Therefore, one seeks to minimize $F(s,\theta(n))= \ee_X\lbr f(s, X_n(s,\theta(n)))\rbr$ (cf.~\cite{yan1992stochastic}, \cite[Chapter 2.5]{KY03}).

\subsection{Main Results}

The aim is to solve the regime-switching discrete stochastic optimization problem~(\ref{eq:SOP}). Inspired by fictitious play learning rules in game theory~\cite{fudenberg1998theory}, we propose a
class of adaptive search algorithms that distributes the search and evaluation functionalities efficiently. The proposed scheme can be described as follows: At each iteration $n$, a state $s(n)$ is sampled from the search space $\M$. The sample $s(n)$ is taken according to a randomized strategy, i.e., a probability distribution on the set $\M$, that minimizes some perturbed variant of the expected objective function based on the beliefs developed thus far. This randomized strategy is referred to \emph{smooth best-response sampling strategy}. The perturbation term in fact simulates the search or exploration functionality essential in learning
the expected stochastic behavior
at various states. The objective function is then simulated at the sampled state $f_n(s(n),\theta(n))$. Finally, the simulation data is fed into a constant step-size stochastic approximation algorithm to update beliefs.

The convergence analysis in Theorem~\ref{theorem:Main} proves that if the underlying hypermodel $\lbr \theta(n)\rbr$ evolves on the same timescale as the
the proposed adaptive search scheme, the most frequently visited state
tracks the set of global optima.
Put differently, the algorithm spends most of its effort simulating the system at the global optima. This is desirable since, in many practical applications, the system has to
be operated in the sampled configuration to measure performance. It is further shown that the proportion of time spent in non-optimal states is inversely proportional to how far their objective function values are from the global minima.
The proposed algorithm relies only on the simulation data and does not require detailed information about the system model, hence, can be used directly as an on-line controller.
The proposed algorithm can, as well, be deployed in static discrete stochastic optimization problems (i.e., when $\theta(n)$ is fixed); see Sec.~\ref{sec:static-DSO} for the related discussion.

The main features of this work are:
\paragraph*{1) Correlated data} We allow for time correlation in samples $f_n(s,\theta(n))$ that is more realistic, whereas most discrete stochastic optimization algorithms assume that the samples are i.i.d.
\paragraph*{2) Adaptive search} The proposed algorithm
tracks the optima as the underlying parameters in the discrete stochastic optimization problem evolve over time. This is in contrast to most existing algorithms that are designed to locate the optima under static settings.
\paragraph*{3) Matched timescale} It is well known that, if the hypermodel $\theta(n)$ changes too drastically, there is no chance one can track the time-varying optima.
(Such a phenomenon is known as trackability; see~\cite{benveniste1990adaptive} for related discussions.) On the other hand, if $\theta(n)$ evolves on a slower timescale as compared to the adaptive search algorithm, it can be approximated by a constant on the fast timescale, hence, its variation is ignored. In this work, we consider the more difficult case where
$\theta(n)$ evolves on the \emph{same} timescale as the adaptive search algorithm and prove that the proposed scheme properly tracks the time varying optima.

Note that the proposed scheme does not assume a Markovian structure for the time-evolution of the objective function. The Markovian switching assumption is only used in our performance analysis that proceeds as follows: First, by a combined use of weak convergence methods~\cite{KY03} and treatment on Markov switched systems~\cite{yin1998continuous,YinZ10}, Theorem~\ref{theorem:discrete-continuous} in Sec.~\ref{sec:asymptotics} shows that the limit system for the discrete time iterates of the proposed algorithm is a randomly switching ordinary differential equation (ODE) modulated by a continuous time Markov chain. (This is in contrast to the standard treatment of stochastic approximation algorithms, where the limiting dynamics converge to a deterministic ODE.) By using multiple Lyapunov function methods for randomly switched systems \cite{chatterjee2007stability,chatterjee2007switch}, Theorem~\ref{theorem:Stability} in Sec.~\ref{sec:stability} proves that the limit switching ODE is asymptotically stable almost surely. Finally, Sec.~\ref{sec:final_convergence} shows that tracking the global attractors set of the derived limit system provides the necessary and sufficient condition to conclude both tracking and efficiency properties of the adaptive search algorithm.

\subsection{Literature}
This work is closely connected to the literature on random search methods;
see~\cite{andradottir2006overview} for a discussion.
Some random search methods spend significant effort to simulate each newly visited state at the initial stages to obtain an
estimate of the objective function. Then, deploying a deterministic optimization mechanism, they search for the global optimum; see~\cite{rubinstein1993discrete,chen2001stochastic,kleywegt2002sample,homem2003variable}. The adaptive search algorithm in this paper is related to another class, namely, discrete stochastic approximation methods~\cite{KY03,spall2003introduction}, which distribute the simulation effort through time, and proceed cautiously based on the limited information available at each time. Algorithms from this class primarily differ in the choice of the sampling strategy.
Examples of sampling strategies can be classified as : i) point-based, leading to methods such as simulated annealing~\cite{gutjahr1996simulated,prudius2012averaging}, tabu search~\cite{glover1997tabu}, stochastic ruler~\cite{alrefaei2005discrete}, stochastic comparison and descent algorithms~\cite{andradottir1996global,andradottir1999accelerating,YKI04,andradottir2009balanced}, ii) set-based, leading to methods such as branch-and-bound~\cite{norkin1998branch}, nested partitions~\cite{shi2000nested}, stochastic comparison and descent algorithms~\cite{hong2006discrete}, and iii) population-based, leading to methods such as genetic algorithms.

Another related body of research pertains to the multi-armed bandit problem~\cite{auer2002finite}, which is concerned with optimizing the cumulative objective function values realized over a period of time, and the pure exploration problem~\cite{audibert2010best}, which involves finding the best arm after a given number of arm pulls.
These methods seek to minimize some regret measure and, similar to the random search methods, usually assume that the problem is static in the sense that the arms' reward distributions are fixed over time\footnote{See~\cite{nonstat-band-2} for upper confidence bound policies for non-stationary bandit problems.}. Further, empirical numerical studies in Sec.~\ref{sec:numerical-example} reveal that bandit-based algorithms such as upper confidence bound (UCB)~\cite{auer2002finite} exhibit reasonable efficiency
only when the size of the search space is relatively small.

\subsection{Organization}
The rest of the paper is organized as follows: Sec.~\ref{sec:DSO-prob} formalizes the main assumptions posed on the problem. In Sec.~\ref{sec:distributed-setting}, the adaptive search scheme is presented and the main theorem of the paper entailing the tracking and efficiency properties is given. Sec.~\ref{sec:proof_main} gives the proof of the main theorem. Finally, numerical examples are provided in Sec.~\ref{sec:numerical-example} followed by the concluding remarks in Sec.~\ref{sec:conclusion}. The proofs are relegated to the Appendix for clarity of presentation.

\section{Main Assumptions}
\label{sec:DSO-prob}
This section formalizes the main assumptions posed on the regime-switching discrete stochastic optimization problem~(\ref{eq:SOP}):

\paragraph{Hypermodel $\theta(n)$} A typical method for analyzing the performance of an adaptive algorithm is to postulate a hypermodel for the underlying time variations~\cite{benveniste1990adaptive}. Here, we assume that all time-varying underlying parameters in the problem are finite-state and absorbed to a vector, indexed by $\theta\in\mathcal{Q}$, whose dynamics follow a discrete-time Markov chain with infrequent jumps. Condition~(A1) below formally characterizes the hypermodel.

\vspace{0.1cm}
\begin{itemize}
\setlength{\itemindent}{1em}
    \item[(A1)] Let $\lbr \theta(n)\rbr$ be a discrete-time Markov chain with finite state space $\mathcal{Q} = \lbr 1,2,\ldots,\Theta\rbr$ and transition probability matrix\footnote{We assume that the initial distribution of the hypermodel $\boldsymbol{p}_0 = [p_{_{0,i}}]_{i\in\mathcal{Q}}$, where  $P(\theta(0) = i) = p_{_{0,i}} \geq 0$ and $\boldsymbol{p}_0 \one_{\Theta} = 1$, is independent of $\vare$.}
\beq{eq:Markov-transition}
P^{\vare} := I +\vare Q.
\eeq
Here, $\vare>0$ is a small parameter, $I$ denotes the $\Theta\times\Theta$ identity matrix, and $Q = \lb q_{ij}\rb \in \rr^{\Theta\times\Theta}$ is the generator of a continuous-time Markov chain satisfying
\begin{equation}
\label{eq:Markov_cont}
q_{ij}\geq0 \;\mathrm{for}\; i\neq j, \;
|q_{ij}|\leq 1\;\forall i,j\in\mathcal{Q}\footnotemark,\; Q\one_{\Theta} = \mathbf{0},
\end{equation}
\footnotetext[3]{This is without loss of generality. Given an arbitrary non-absorbing Markov chain with transition probability matrix $P^{\rho} = I+\rho\hat{Q}$, one can form $Q = \frac{1}{q_{\textmd{max}}}\cdot\hat{Q}$, where $q_{\textmd{max}} = \max_{i\in\M} |\bar{q}_{ii}|$. To ensure that the two Markov chains generated by $Q$ and $\bar{Q}$ evolve on the same timescale, $\varepsilon = \rho\cdot q_{\textmd{max}}$.}

\vspace{-0.3cm}
\noindent
where $\one_{\Theta} = \lb 1,\ldots,1\rb_{\Theta\times 1}$ and $Q$ is irreducible.
\end{itemize}
\vspace{0.1cm}
Choosing $\vare$ small enough ensures that the entries of $P^{\vare}$ in (\ref{eq:Markov-transition}) are non-negative.
The use of the generator $Q$ also makes the row sum of $P^{\vare}$ be one. Due to the dominating identity matrix in (\ref{eq:Markov-transition}), $\lbr \theta(n)\rbr$ varies slowly with time.

%
\vspace{0.1cm}

\paragraph{Simulation Data $f_n\left(s,\theta(n)\right)$}
Let $\ee_\ell$ denotes the conditional expectation given $\mathcal{F}_\ell$, the $\sigma$-algebra generated by $\lbr f_n(s,\theta(n)), s\in\M, \theta(n) : n < \ell\rbr$. We make the following assumptions.

\vspace{0.1cm}
\begin{itemize}
\setlength{\itemindent}{1em}
    \item[(A2)] For each $s\in\M$ and $\theta\in\mathcal{Q}$, $\lbr f_n\left(s,\theta\right)\rbr$ is a sequence of bounded real-valued random variables. Moreover, for any $\ell\ge 0$,
        %
        \begin{equation}
        \label{eq:A2}
        \hspace{-0.1cm}
        \frac{1}{n} \sum^{n+\ell-1}_{\tau=\ell}
        \ee_\ell f_\tau\left(i,j\right) \to F\left(i,j\right) \
        \hbox{ in probability as}\ n\to\infty 
        ,
        \end{equation}
        for all $i\in\M$ and $j\in\mathcal{Q}$, where
        $F\left(s,\theta\right) = \ee \lbr f_n(s,\theta)\rbr$; see Sec.~\ref{sec:intro-problem}.
        %
\end{itemize}
\vspace{0.1cm}
The above condition allows us to work with correlated processes whose remote past and distant future are asymptotically independent. Examples include: the sequence of i.i.d. random variables with (asymptotically) uniformly bounded variance, or a class of random variables (not necessarily i.i.d.) that satisfy the large deviations principle (cf.~\cite{hong2006discrete,yakowitz2000global}), e.g., moving average and stationary auto-regressive processes.

Finally, we impose the following condition on the hypermodel $\theta(n)$: Let $\mu$ denote the adaptation rate of the adaptive search algorithm; see~(\ref{eq:avg-normalized}) or~(\ref{eq:regret-update}). Then,
\vspace{0.1cm}
\begin{itemize}
\setlength{\itemindent}{1em}
    \item[(A3)] $\vare = \mu$ in the transition probability matrix $P^\vare$.
\end{itemize}
\vspace{0.1cm}
Condition (A3) states that time variations of the parameters underlying the discrete stochastic optimization problem (\ref{eq:SOP}) occur at the same timescale as the updates in the proposed adaptive search algorithm.

\vspace{0.1cm}
\begin{Remark}
It is important to stress that the hypermodel $\theta(n)$ is not used in the algorithm proposed in this paper. The algorithm does not require knowledge of $\theta(n)$ or its parameters. The hypermodel is used only in the analysis of the algorithm. We are interested in determining if the algorithm can track time-varying optima that evolve according to a slow Markov chain. Since
$\theta(n)$ is unobservable, we suppress the dependence of $f_n(s,\theta(n))$ on
it and, with slight abuse of notation, denote it by $f_n(s)$ in what follows.
\end{Remark}

\section{Tracking the Global Optima: Algorithm and Main Results}
\label{sec:distributed-setting}
In this section, we introduce a stochastic approximation algorithm that, relying on smooth best-response strategies~\cite{hofbauer2002global,BHS06}, prescribes how to sample from the search space so as to efficiently learn and track the evolving set of global optima $\go(\theta(n))$. To this end, we define the smooth best-response procedure based on consecutive observations $\lbrace f_n(s_n) \rbrace_{n\geq 0}$ and outline its distinct properties in Sec.~\ref{sec:smooth-bestresponse}. We then present the proposed adaptive discrete stochastic optimization algorithm in Sec.~\ref{sec:algorithm} followed by the main result of the paper that shows, if one employs the proposed algorithm and the time-varying underlying parameters evolve on the same timescale as the
the stochastic approximation algorithm, the algorithm efficiently tracks the set of global optima $\go(\theta(n))$.

\subsection{Smooth Best-Response Sampling Strategy}
\label{sec:smooth-bestresponse}
Consider a learning scenario where one repeatedly samples from the search space, denoted by $s(n)\in\M$, at discrete times $n = 1,2,\ldots$ and obtains $f_n(s(n))$ via simulation or measurement. We postulate that $s(n)$ is chosen according to a randomized sampling strategy $\p(n) = \left( p_1(n),\ldots,p_S(n)\right)$ that belongs to the simplex of probability distributions over the search space
\begin{equation}
\label{eq:delta-M}
\Delta\M = \lbr \p\in\RR^{S}; p_i\geq 0, \sum_{s\in\M} p_i = 1\rbr.
\end{equation}
Based only on the collected observations $\lbrace f_\tau(s(\tau)) : \tau\leq n\rbrace$ up to time~$n$, define the vector of weighted average objective function values $\f(n) = \lb\tilde{f}_1(n),\ldots,\tilde{f}_S(n)\rb'\in\RR^S$, where $v'$ denotes the transpose of $v$, and
\begin{equation}
\label{eq:avg-normalized}
\tilde{f}_i(n) = \mu \sum_{\tau\leq n} (1-\mu)^{n-\tau} \frac{f_\tau(s(\tau))}{p_i(\tau)}\cdot I_{\lbr s(\tau) = i\rbr}, \quad\forall i\in\M.
\end{equation}
In~(\ref{eq:avg-normalized}), $I_{\lbrace \cdot\rbrace}$ denotes the indicator function, and the normalization factor $1/p_i(\tau)$ makes the length of the periods that each states $i$ is chosen comparable to other states. The discount factor $\mu$ places more weight on recent observations and is necessary as the algorithm is deemed to track time-varying minima. Note further that (\ref{eq:avg-normalized}) only relies on the actual measurements or simulation data $\tilde{f}_\tau(s(\tau))$ recorded (e.g. from the system performance) and does not require the system model nor the realizations of $\lbr\theta_{\tau}\rbr$.
The smooth best-response sampling strategy is then defined as follows.
\vspace{0.1cm}

\begin{Definition}
\label{def:smooth-best}
Choose a function $\rho(\boldsymbol{\sigma}):\textmd{int}(\Delta \M)\to \RR$, where $\textmd{int}(G)$ denotes the interior
 of $G$
 and $\Delta \M$ is defined in~(\ref{eq:delta-M}), such that
\begin{itemize}
    \item[i)] $\rho\cd$ is $\mathcal{C}^1$ (i.e., continuously differentiable), strictly concave, and $|\rho|\leq 1$;
    \item[ii)] $\|\nabla\rho(\boldsymbol{\sigma})\|\to\infty$ as $\boldsymbol{\sigma}$ approaches the boundary of $\Delta \M$, i.e.,
    $$\lim_{\boldsymbol{\sigma}\to\partial \left(\Delta\M\right)} \left\| \nabla\rho(\boldsymbol{\sigma})\right\| = \infty,$$
    where $\|\cdot\|$ denotes the Euclidean norm, and $\partial(\Delta\M)$ represents the boundary of simplex $\Delta \M$.
\end{itemize}
The \emph{smooth best-response sampling} strategy is then given by
\begin{equation}
\label{eq:smooth-br}
\br\big(\f\big) := \operatorname*{arg\,min}_{\boldsymbol{\sigma} \in\Delta \M} \sum_{i\in\M} \sigma_i \tilde{f}_i - \gamma \rho(\boldsymbol{\sigma}),\quad 0<\gamma<\hat{\gamma}.
\end{equation}
\end{Definition}

The conditions imposed on the perturbation function $\rho\cd$ leads to the following distinct properties of the resulting strategy:
\begin{itemize}
    \item[i)] The strict concavity condition ensues the uniqueness of $\br\big(\f\big)$;
    \item[ii)] The boundary condition implies $\br\big(\f\big)$ belongs to the interior of the simplex $\Delta \M$.
\end{itemize}
The smooth best-response strategy is inspired by leaning algorithms in games~\cite{hofbauer2002global,fudenberg1998theory}.
It exhibits exploration using the idea of adding a random value to the belief about the objective function values associated with each state. (This is in contrast to picking states at random with a small probability, as is common in game-theoretic learning and multi-armed bandit algorithms.) Such exploration is natural in any learning scenario. The results of~\cite[Theorem~2.1]{hofbauer2002global} show that, regardless of the distribution of the random values, a deterministic representation of the form (\ref{eq:smooth-br}) can be obtained for the pure best-response strategy resulted from adding random values to the beliefs $\f(n)$. Further,
the smooth best-response strategy constructs a genuine randomized strategy. This is an appealing feature since it circumvents the discontinuity inherent in algorithms of pure best-response type (i.e., $\operatorname*{arg\,max}_{i \in \M}  \tilde{f}_i$), where small changes in the beliefs $\f(n)$ can lead to an abrupt change in the behavior of the algorithm. Such switching behavior in the dynamics of the algorithm
complicates the convergence analysis.
\vspace{0.1cm}
\begin{Remark}
\label{remark-3-1}
An example of the function $\rho\left(\cdot\right)$ in Definition~\ref{def:smooth-best} is the \emph{entropy function}~\cite{fudenberg1998theory,fudenberg1999conditional}
\begin{displaymath}
\rho\left(\boldsymbol{\sigma}\right) = -\sum_{i\in\M} \sigma_i \ln\left(\sigma_i\right),
\end{displaymath}
which gives rise to the smooth best-response strategy
\begin{equation}
\bbr_i\big(\tilde{\boldsymbol{f}}\big) = \frac{\exp\big(-\tilde{f}_i/\gamma\big)}{\sum_{j\in\M}\exp\big(-\tilde{f}_j/\gamma\big)}.
\end{equation}
Such a strategy is also used in the context of learning in games, widely known as logistic fictitious-play~\cite{fudenberg1995consistency} or logit choice function~\cite{hofbauer2002global}.
\end{Remark}

\subsection{Adaptive Discrete Stochastic Optimization Algorithm}
\label{sec:algorithm}
We now proceed to present the stochastic approximation algorithm proposed for tracking the set of global optima $\go(\theta(n))$. The adaptive discrete stochastic optimization algorithm can be simply described as an adaptive sampling scheme. Relying on the beliefs developed about the objective function values at each states, it prescribes how to sample from the search space $\M$ so as to efficiently (in terms of the amount of effort spent on simulating non-promising states) track the global optima $\go(\theta(n))$. We then
simulate $f_n(s(n))$
at the sampled state $s(n)$ and use a stochastic approximation algorithm to update beliefs $\f(n)$ and, accordingly, the sampling strategy. The proposed algorithm relies only on the simulation
data and is efficient in the sense that it requires minimum effort per iteration---it needs only one simulation, as compared to, e.g.,
two in~\cite{andradottir1996global}. Yet, as evidenced by the numerical example in Sec.~\ref{sec:numerical-example}, it guarantees performance gains in terms of tracking speed.

The adaptive discrete stochastic optimization algorithm is summarized below:

\vspace{0.1cm}
\noindent
\textbf{Algorithm~1}:

\noindent
\emph{Aim.} Generate a sequence $\lbr s(n)\rbr$ that provides an estimate of the time-varying global optima.

\vspace{0.2cm}
\textbf{\textit{Step 0)} Initialization:} Choose $\rho\cd$ to satisfy the conditions of Definition~\ref{def:smooth-best} and set the exploration parameter $\gamma>0$. \newline Initialize $\f(0) = \mathbf{0}_{S}$.

\textbf{\textit{Step 1)} State Sampling:} Select state $s(n) \sim \br\big(\f(n)\big)$; see (\ref{eq:smooth-br}).

\textbf{\textit{Step 2)} State Evaluation:} Simulate or measure $f_n(s(n))$.

\textbf{\textit{Step 3)} Belief Update:} Update the $S$-dimensional vector
\beq{eq:regret-update}
\f(n+1) = \f(n) + \mu\left[\boldsymbol{g}\left(s(n),\f(n)\right) - \f(n)\right],
\eeq
where $\boldsymbol{g}\big(s(n),\f(n)\big)$ is a column vector with elements
\beq{eq:A-elements}
g_{i}\left(s(n),\f(n)\right) := \frac{f_n(s(n))}{\bbr_i\big(\tilde{\boldsymbol{f}}(n)\big)} \cdot I_{\lbr s(n) = i\rbr}.
\eeq

\textbf{\textit{Step 4)} Recursion:} Set $n\leftarrow n+1$ and go to Step 1.

\vspace{0.3cm}
\begin{Remark}
1)
Note that
the dynamics of $\lbr \theta(n)\rbr$
 do not enter implementation of the algorithm, and is only used in the tracking analysis in Sec.~\ref{sec:asymptotics}. In particular, Theorem~\ref{theorem:Main} shows that Algorithm~1 can successfully track the time-varying optima if they vary according to the hymeromodel $\theta(n)$.

2) If $\lbr \theta(n)\rbr$ was observed, one could form and update $\f_\theta(n)$ independently for each $\theta\in\mathcal{Q}$, and use $\br\big(\f_{\theta^\prime}(n)\big)$ to select $s(n)$ once the system switched to $\theta^\prime$. It can then be shown that the sequence $\lbrace s(n)\rbrace$ spends most of its time in the global minima, irrespective of the switching, for all $\theta\in\mathcal{Q}$.

3)
Larger values of $\gamma$ increase the exploration weight versus exploitation, hence, decreases the time spent in $\go(\theta(n))$.
\end{Remark}

\subsection{Main Result: Tracking the Regime-Switching Minima Set}
\label{sec:main-results}

To analyze the tracking capability of the above adaptive discrete stochastic optimization algorithm, define two diagnostics that will be used subsequently:

(i) \emph{Regret $r(n)$}:
\begin{equation}
\label{eq:regret-org}
r(n) := \bar{f}(n) - F_{\min}(\theta(n)),
\end{equation}
where
\begin{align}
\label{eq:average_1}\bar{f}(n) &:= \mu\sum_{\tau\leq n} (1-\mu)^{n-\tau} f_\tau(s(\tau)),\\
\label{eq:average_2}F_{\min}(\theta) & := \min_{s\in\M} F(s,\theta).
\end{align}
%
Here, $\lbr s(n)\rbr$ is the sequence of states prescribed by the discrete stochastic optimization algorithm, and $\bar{f}(n)$ represents the expected realized objective function value up to sampling period $n$. Thus, the regret $r(n)$ quantifies the tracking capability of the algorithm.


(ii) \emph{Empirical Sampling Distribution:} To study efficiency of the adaptive discrete stochastic optimization algorithm, we define the empirical sampling distribution vector $\z(n)\in\RR^S$ as
\beq{eq:empirical_frequency_local}
\z(n) := \mu\sum_{\tau\leq n}\left(1-\mu\right)^{n-\tau}\boldsymbol{e}_{s(\tau)},
\eeq
where $\boldsymbol{e}_i\in\RR^S$ denotes the unit vector with the $i$th element being equal to one. Therefore, $z_i(n)$ records the percentage of iterations that state $i$ was sampled and simulated up to time $n$. Efficiency of a discrete stochastic optimization algorithm is defined as the percentage of time that states within the set of global optima
are sampled. For each $\bar{\theta}\in\mathcal{Q}$, the efficiency is thus quantified by $\sum_{i\in\go(\bar{\theta})} z_i(n)$. In~(\ref{eq:empirical_frequency_local}), $\mu$ serves as the forgetting factor to facilitate adaptivity to the evolution of underlying parameters.

Before proceeding with the main theorem, define the continuous time interpolated sequence of iterates
\begin{equation}
\label{eq:interpol-z-r}
\begin{array}{c}
\z^{\mu}(t) = \z(n)\\
r^{\mu}(t) = r(n)
\end{array}
\quad \textmd{for} \quad t\in\left[n\mu,(n+1)\mu\right],
\end{equation}
and let
\beq{eq:simplex-epsilon}
\Upsilon^\eta(\theta) = \lbr \boldsymbol{\pi}\in \Delta\M ; \sum_{i\in\M}\pi_i\lb f(i,\theta) - F_{\min }(\theta)\rb\leq \eta\rbr.
\eeq
The following theorem
asserts that the sequence $\lbr s(n)\rbr$ generated by Algorithm~1 tracks the regime-switching minima set $\go(\theta(n))$ and spends most of its effort on simulating
$\go(\theta(n))$. In what follows, $\Rightarrow$ denotes weak convergence.\footnote{Weak convergence is a generalization of convergence in distribution to a function space~\cite{KY03}; see also Sec.~\ref{sec:asymptotics} of this paper for a precise definition.} Note that when a sequence converges weakly to a constant, it also converges in probability to that constant.
\vspace{0.1cm}
\begin{Theorem}
\label{theorem:Main}
Suppose (A1), (A2), and (A3) hold. Let $q(\mu)$ be any sequence of real numbers satisfying $q(\mu)\rightarrow\infty$ as $\mu \rightarrow 0$. Then, for any $\eta>0$, there exists $\bar{\gamma} > 0$ such that, if $\gamma\leq\bar{\gamma}$ in~(\ref{eq:smooth-br}), as $\mu\to 0$\footnote{We assume the initial values $\z(0)$ and $r(0)$ are independent of the step-size $\mu$ for simplicity. Otherwise, if $\z(0) = \z^{\mu}(0)$ and $r(0) = r^{\mu}(0)$, we require that $\z^{\mu}(0)$ and $r^{\mu}(0)$ converge weakly to $\z(0)$ and $r(0)$, respectively.}:
\begin{enumerate}
\item \emph{Tracking:} $(
r^\mu(\cdot+q(\mu)) - \eta
)^+\Rightarrow 0$, where $
x^+ = \max\lbr 0,x\rbr$.
\item \emph{Efficiency:} $\z^{\mu}(\cdot+q(\mu))\Rightarrow\Upsilon(\theta\cd)$ in the sense that:
\beq{eq:main-thrm-1}
\begin{split}
&d\big(\z^{\mu}(\cdot + q(\mu)),\Upsilon^\eta(\theta\cd)\big)\\
&\hspace{1cm}=\inf_{\boldsymbol{\pi}\cd\in\Upsilon^\eta(\theta\cd)} \big|\z^{\mu}(\cdot + q(\mu))-\boldsymbol{\pi}\cd\big| \Rightarrow 0,
\end{split}
\eeq
where $d(\cdot,\cdot)$ is the usual distance function, and $\theta(\cdot)$ is a continuous time Markov chain with generator $Q$; see~(A1).
\end{enumerate}

\begin{IEEEproof}
The
proof
uses martingale averaging techniques to show that the limit behavior converges weakly to a switched Markovian ordinary differential equation~(ODE). Then, stability of the switched ODE is established and the global attractor set is shown to represent the global minina set. The detailed proof is in Sec.~\ref{sec:proof_main}.
\end{IEEEproof}
\end{Theorem}

\paragraph*{Interpretation of Theorem~\ref{theorem:Main}}
The above theorem addresses both \emph{tracking capability} and \emph{efficiency} of Algorithm~1: Part 1) evidences both \emph{consistency} and \emph{attraction} to the set $\go(\theta\cd)$ by looking at the continuous time interpolation of worst case regret $r^\mu\cd$. It shows that $r^\mu(t)$ stays infinitely often less that $\eta$ as $\mu\to 0$ and $t\to\infty$. (This result is similar to the \emph{Hannan consistency} notion~\cite{hannan1957approximation} in repeated games, however, in a regime-switching setting.) Part~2) concerns efficiency by showing that the algorithm eventually spends most of its effort on simulating $\go(\theta\cd)$ and adapts to its time variations. In particular, the proportion of time spent simulating states $s\notin\go(\theta\cd)$ is inversely proportional to how far their objective value is from the global minimum. Note that Part 2) claims convergence to a set, rather than a point in the set.

\vspace{0.1cm}

The following corollary is a direct consequence of Theorem~\ref{theorem:Main}. It asserts that the continuous time interpolation of the most frequently  visited state converges weakly to the set of global minima.
\vspace{0.1cm}
\begin{Corollary}
\label{corrolary}
Denote the most frequently visited state by
\begin{equation}
\label{eq:max}
s_{_{\textmd{max}}}(n) = \argmax_{i\in\M}\ z_i(n),\nonumber
\end{equation}
where $\bar{z}_i(n)$ is the $i$th component of $\z(n)$ defined in~(\ref{eq:empirical_frequency_local}). Define the continuous time interpolated sequence
\begin{equation}
s_{_{\textmd{max}}}^{\mu}(t) = s_{_{\textmd{max}}}(n)\quad \textmd{for}\quad t\in\left[n\mu,(n+1)\mu\right].\nonumber
\end{equation}
Then, under (A1)--(A3) and the conditions of Theorem~\ref{theorem:Main}, $s_{_{\textmd{max}}}^{\mu}(\cdot + q(\mu))$ converges weakly to the regime-switching minima set $\go(\theta\cd)$ as $\mu\to\infty$.
\end{Corollary}
\vspace{0.1cm}

Note that, to foster adaptivity to the time variations of the hypermodel $\lbr \theta(n)\rbr$, Algorithm~1 selects non-optimal states with some small probability. Thus, one would not expect $\lbr s(n)\rbr$ to converge to $\go(\theta(n))$. In fact, $\lbr s(n)\rbr$ may visit each element of $\M$ infinitely often. Instead, the strategy implemented by following Algorithm~1 ensures the empirical frequency of sampling from $\M\backslash\go(\theta(n))$ stays very low.

\subsection{Static Discrete Stochastic Optimization}
\label{sec:static-DSO}

Suppose $\theta(n) = \bar{\theta}$ is fixed in~(\ref{eq:SOP}). The discrete stochastic optimization problem then reduces to
\beq{eq:DSO-static}
\min_{s\in\M} F\left(s\right) = \ee\lbr f_n\left(s,\bar{\theta}\right)\rbr,\nonumber
\eeq
and is \emph{static} in the sense that the set $\go(\bar{\theta})$ of global minima does not evolve with time. Although not being the focus of this paper, one can use the results of~\cite{benaim2012consistency} to show that if the exploration factor $\gamma$ in (\ref{eq:smooth-br}) decreases to zero sufficiently slowly, the sequence $\lbr s(n)\rbr$ converges almost surely to $\go(\bar{\theta})$.

More precisely, consider the following modifications to Algorithm~1:
\begin{enumerate}
    \item[(i)] The constant step-size $\mu$ in (\ref{eq:regret-update}) is replaced by decreasing step-size $\mu_n = \frac{1}{n+1}$;
    \item[(ii)] The exploration factor $\gamma$ in (\ref{eq:smooth-br}) is replaced by $\frac{1}{n^\alpha}$, where $0<\alpha<1$.
\end{enumerate}
Define the sequence of interpolated process $s^{n}(t)$, $n=0,1,\ldots$:
\begin{equation}
\begin{split}
&\hspace{0.55cm}s^{0}(t) = s(n)\ \textmd{for} \ t\in[t_n,t_{n+1}),\\
&s^{n}(t) = s^{0}(t+t_n)\ \textmd{for} \ -\infty<t<\infty,
\end{split}\nonumber
\end{equation}
where $t_n = \sum_{\tau=0}^{n-1}\mu_\tau$. Let
$q(n)$ be any sequence of real numbers satisfying $q(n)\to\infty$ as $n\to\infty$. Then, if $\lbr s(n)\rbr$ is chosen according to Algorithm~1, $s^{n}(\cdot+q(n)) \xrightarrow{\textmd{a.s.}} \go(\bar{\theta})$ as $n\to\infty$ in the sense that $d\big(s^{n}(\cdot+q(n)),\go(\bar{\theta})\big)\xrightarrow{\textmd{a.s.}} 0$.

By the above construction, the sequence $\lbr s(n)\rbr$ will eventually become reducible with singleton communicating class $\go(\bar{\theta})$. That is, $\lbr s(n)\rbr$ eventually spends all its time in $\go(\bar{\theta})$. This is in contrast with Algorithm~1 in the regime-switching setting.

\section{Proof of Theorem~\ref{theorem:Main}: Tracking regime-switching Global Minima}
\label{sec:proof_main}
This section presents the proof of the main result and is organized into four subsections: We start by showing in Sec.~\ref{sec:asymptotics} that the limit system associated with the discrete time iterates $\big(\f(n),r(n)\big)$ is a Markovian switching system of interconnected ODEs. Next, Sec.~\ref{sec:stability} proves that such a limit system is globally asymptotically stable with probability one and characterizes its global attractors. Accordingly, we conclude asymptotic stability of the interpolated process associated with $\big(\f(n),r(n)\big)$ in Sec.~\ref{sec:asymp_stability}, and prove that the the discrete time iterates mimicking such limit dynamics is attracted to the set of global minima. Finally, Sec.~\ref{sec:final_convergence} uses the results obtained thus far to conclude efficiency of Algorithm~1.

\subsection{Weak Convergence to Markovian Switching ODE}
\label{sec:asymptotics}
In this subsection, we use weak convergence methods to derive the limit dynamical system associated with the iterates $\big(\f(n),r(n)\big)$. Before proceeding further, let us recall some definitions and notation:

Let $Z(n)$ and $Z$ be $\RR^{\mathfrak{s}}$-valued random vectors. We say $Z(n)$ converges weakly to $Z$ ($Z(n) \Rightarrow Z$) if for any bounded and continuous function $\psi(\cdot)$, $E\psi(Z(n))\to E\psi(Z)$ as $n\to \infty$. We also say that the sequence $\{Z(n)\}$ is tight if for each $\eta>0$, there exists a compact set $K_\eta$ such that $P(Z(n)\in K_\eta)\ge 1-\eta$ for all $n$. The definitions of weak convergence and tightness extend to random elements in more general metric spaces. On a complete separable metric space,  tightness is equivalent to relative compactness, which is known as Prohorov's Theorem~\cite{billingsley1968convergence}. By virtue of this theorem, we can extract convergent subsequences when tightness is verified. In what follows, we use a martingale problem formulation to establish the desired weak convergence. To this end, we first prove tightness. The limit process is then characterized using a certain operator related to the limit martingale problem. We refer the reader to \cite[Chapter 7]{KY03} for further details on weak convergence and related matters.

\def\fhat{\hat{\boldsymbol{f}}}
\def\bf{\boldsymbol{F}}
Define
\begin{equation}
\label{eq:F-theta}
\boldsymbol{F}(\theta) = \lb F(1,\theta),\cdots,F(S,\theta)\rb',
\end{equation}
and let
\begin{equation}
\label{eq:f-hat}
\fhat(n) := \f(n) - \bf(\theta(n))
\end{equation}
denote the deviation error in tracking the true objective function values via the simulation data at time $n$. Let further
\begin{equation}
\label{eq:X_def}
\boldsymbol{X}(n) := \left[
\begin{array}{c}
\fhat(n)\\
r(n)
\end{array}
\right].
\end{equation}
It can be easily verified that $\boldsymbol{X}(n)$ satisfies the recursion
\begin{equation}
\label{eq:X-recursion}
\begin{split}
\boldsymbol{X}(n+1) &= \boldsymbol{X}(n) + \mu\lb \boldsymbol{A}_n(s(n)) - \boldsymbol{X}(n)\rb \\
&+ \left[
\begin{array}{c}
\boldsymbol{F}(\theta(n)) - \boldsymbol{F}(\theta(n+1))\\
F_{\min}(\theta(n)) - F_{\min}(\theta(n+1))
\end{array}
\right],
\end{split}
\end{equation}
where
\begin{equation}
\label{eq:A1-def}
\begin{split}
&\hspace{0.7cm}\boldsymbol{A}_n(s(n)) = \left[
\begin{array}{c}
\hat{\boldsymbol{g}}\left(s(n),\;\fhat(n)\right) - \bf(\theta(n))\\
f_n(s(n)) - F_{\min}(\theta(n))
\end{array}
\right],\\
&\hat{\boldsymbol{g}} = [\hat{g}_1,\cdots,\hat{g}_S]',\; \hat{g}_i = \frac{f_n(s(n))}{b_i^\gamma\left(\fhat(n) + \bf(\theta(n))\right)}\cdot I_{\lbr s(n) = i\rbr},
\end{split}
\end{equation}
and $F_{\min}\cd$ and $\boldsymbol{F}\cd$ are defined in~(\ref{eq:average_2}) and~(\ref{eq:F-theta}), respectively. As is widely used in the analysis of stochastic approximations, we consider the piecewise constant continuous time interpolated processes
\begin{equation}
\label{eq:interpolations}
\barray
\boldsymbol{X}^\mu(t) = \boldsymbol{X}(n),\\
\theta^\mu(t) = \theta(n),
\earray
\quad \textmd{for} \;\; t\in[n\mu,(n+1)\mu).
\end{equation}
In what follows, we use $D\big([0,\infty):
\tilde G \big)$ to denote the space of functions that are defined in $[0,\infty)$ taking values in
$\tilde G$,
and are right continuous and have left limits with Skorohod topology (see \cite[p. 228]{KY03}). The following theorem characterizes the limit process of the stochastic approximation iterates as a Markovian switching ODE.

\vspace{0.1cm}
\begin{Theorem}
\label{theorem:discrete-continuous}
Consider the recursion~(\ref{eq:X-recursion}) and suppose (A1), (A2), and (A3) hold. As $\mu\rightarrow 0$, the interpolated process $(\boldsymbol{X}^\mu\cd,\theta^{\mu}\cd)$ is tight in $D([0,\infty): \RR^{S+1}\times\mathcal{Q})$ and converges weakly to $(\boldsymbol{X}\cd,\theta\cd)$ that is a solution of the Markovian switched ODE
\begin{equation}
\label{eq:ODE-global}
\frac{d\boldsymbol{X}}{dt} = \boldsymbol{G}(\boldsymbol{X},\theta(t)) - \boldsymbol{X},
\end{equation}
where
\begin{equation}
\label{eq:G}
\boldsymbol{G}(\boldsymbol{X},\theta(t)) = \left[
\begin{array}{c}
\mathbf{0}_S\\
\br\left(\fhat + \bf(\theta(t))\right)\cdot \boldsymbol{F}(\theta(t)) - F_{\min}(\theta(t))
\end{array}
\right].\nonumber
\end{equation}
Here, $\mathbf{0}_S$ denotes an $S\times 1$ zero vector, $\bf(\cdot)$ and $F_{\min}(\cdot)$ are defined in~(\ref{eq:average_2}) and~(\ref{eq:F-theta}), respectively, and $\theta(t)$ denotes a continuous time Markov chain with generator $Q$; see~(A1).
\begin{IEEEproof}
The proof
uses stochastic averaging theory based on~\cite{KY03}; see Appendix~\ref{sec:proof_1} for the detailed argument.
\end{IEEEproof}
\end{Theorem}
\vspace{0.1cm}

The above theorem asserts that the asymptotic behavior of Algorithm~1 can be captured by a dynamical system modulated by a continuous-time Markov chain $\theta(t)$. At any given instance, the Markov chain dictates which regime the system belongs to, and the system then follows the corresponding ODE until the modulating Markov chain jumps into a new state (i.e., the limit system~(\ref{eq:ODE-global}) is only piecewise deterministic).

\vspace{0.1cm}
\begin{Remark}
When $\theta(n)$ evolves on a slower timescale, e.g., $\vare = \mathcal{O}(\mu^2)$ in~(\ref{eq:Markov-transition}), it remains constant in the fast timescale (i.e., the adaptive discrete stochastic optimization algorithm). Therefore, the ODE~(\ref{eq:ODE-global}) will become deterministic.
\end{Remark}

\subsection{Stability Analysis of the Markovian Switching ODE}
\label{sec:stability}
We next proceed to analyze stability and characterize the set of global attractors of the limit system~(\ref{eq:ODE-global}).

Let us start by looking at the evolution of the deviation error $\fhat(t)$ in tracking the objective function values, which forms the first component in any trajectory $\X(t)$ of the limit system. In view of~(\ref{eq:ODE-global})--(\ref{eq:G}), $\fhat(t)$ evolves according to the deterministic ODE
\begin{equation}
\label{eq:f-hat-ode}
\frac{d\fhat}{dt} = -\fhat.\nonumber
\end{equation}
Note that the dynamics of $\fhat(t)$ is independent of the second component of $\X(t)$, namely, the regret $r(t)$. Since the ODE is asymptotically stable, $\fhat(t)$ decays exponentially fast to $\mathbf{0}_S$ as $t\to\infty$.
This essentially establishes that realizing $f_n(s(n))$ provides sufficient information to construct an unbiased estimator of the true objective function values\footnote{It can be shown that the sequence $\big\{\f(n)\big\}$ induces the same asymptotic behavior as the beliefs developed using the brute force scheme~\cite[Chapter 5.3]{pflug1996optimization} about objective function values.}.

Next, substituting the global attractor $\fhat = \mathbf{0}_S$ into the limit switching ODE associated with the regret $r(t)$ (the second component in $\X(t)$), we analyze stability of
\begin{equation}
\label{eq:r-ode}
\frac{dr}{dt} = \br\left(\bf(\theta(t))\right)\cdot \boldsymbol{F}(\theta(t)) - F_{\min}(\theta(t)) - r.
\end{equation}
We start by defining stability of switched dynamical systems; see~\cite[Chapter 9]{YinZ10} and~\cite{chatterjee2007switch} for further details. In what follows, $d(\cdot,\cdot)$ denotes the usual distance function.
\vspace{0.1cm}
\begin{Definition}
\label{def:stability}
Consider the Markovian switched system
\beq{eq:switched system}
\begin{split}
&\hspace{1.9cm}\dot{Y}(t) = f\left(Y(t),\theta(t)\right)\\
&Y(0) = Y_0,\ \theta(0) = \theta_0,\ Y(t)\in\mathbb{R}^r,\ \theta(t)\in\mathcal{Q},
\end{split}
\nonumber
\eeq
where $\theta(t)$ is a continuous time Markov chain with generator $Q$, and
$f(\cdot, i)$ is
locally Lipschitz for each $i\in\mathcal{Q}$. A closed and bounded set $\H \subset \mathbb{R}^n\times \mathcal{Q}$ is:
\begin{enumerate}
    \item \emph{stable in probability} if
    for any $\varrho, \bar\gamma>0$, there is a $\gamma>0$ such that
    $$\mathds{P}\bigg(\sup_{t\geq 0} \ d\big((Y(t), \theta(t)),\H\big)<\bar\gamma \bigg) \ge  1-\varrho,$$
    whenever $d\left((Y_0,\theta_0),\H\right)< \gamma$;
    \item \emph{asymptotically stable in probability} if it is stable in probability and
    $$\mathds{P}\Big(\lim_{t\to \infty} d\big((Y(t),\theta(t)), \H\big)=0\Big)\to 1;$$
   \item \emph{asymptotically stable almost surely} if
   $$\lim_{t\to \infty} d\big((Y(t), \theta(t)), \H\big) =0 \ \hbox{ a.s.}$$
\end{enumerate}
\end{Definition}
\vspace{0.1cm}

Before proceeding with the theorem, let
\beq{eq:R}
\mathbb{R}_{[0,\eta)} = \lbr r\in\RR; 0\leq r < \eta\rbr.
\eeq
We break down the stability analysis of~(\ref{eq:r-ode}) into two steps; First, we examine the stability of each subsystem, i.e., for each $\bar{\theta}\in\mathcal{Q}$ when $\theta(t) = \bar{\theta}$ is fixed. The set of global attractors is shown to comprise $\mathbb{R}_{[0,\eta)}$ for all $\bar\theta\in\mathcal{Q}$. The slow switching condition then allows us to apply the method of multiple Lyapunov functions~\cite[Chapter 3]{liberzon2003switching} to analyze stability of the switched system.
\vspace{0.1cm}
\begin{Theorem}
\label{theorem:Stability}
Consider the limit Markovian switched ODE given in~(\ref{eq:r-ode}). Let $r(0) = r_0$ and $\theta(0) = \theta_0$. For any $\eta>0$, there exists $\bar{\gamma}(\eta)$ such that, if $\gamma<\bar{\gamma}(\eta)$ in~(\ref{eq:smooth-br}), the following results hold:
\begin{enumerate}
    \item If $\theta(t) = \bar{\theta}$ is fixed, the deterministic dynamical system~(\ref{eq:r-ode}) is asymptotically stable., the set $\mathbb{R}_{[0,\eta)}$ is globally asymptotically stable for each $\bar\theta\in\mathcal{Q}$, i.e.,
    \begin{equation}
    \lim_{t\to\infty} d\left(r(t),\mathbb{R}_{[0,\eta)}\right) = 0.
    \end{equation}
    \item For the Markovian switching ODE,
    the set $\mathbb{R}_{[0,\eta)}$ is globally asymptotically stable almost surely.
\end{enumerate}
\end{Theorem}

\begin{IEEEproof}
For detailed proof, see Appendix~\ref{sec:proof_2}.
\end{IEEEproof}

\vspace{0.1cm}

The above theorem states that the set of global attractors of the switching ODE (\ref{eq:r-ode}) is the the same as that for all non-switching ODEs (i.e., when $\theta(t)=\bar\theta\in\mathcal{Q}$ is fixed in (\ref{eq:r-ode})) and constitutes $\mathbb{R}_{[0,\eta)}$. This sets the stage for Sec.~\ref{sec:final_convergence} where attraction to $\mathbb{R}_{[0,\eta)}$ is shown to conclude the desired tracking and efficiency results.

\subsection{Asymptotic Stability of the Interpolated Process}
\label{sec:asymp_stability}

In Theorem~\ref{theorem:discrete-continuous}, we considered $\mu$
small and $n$ large, but $\mu n$ remained bounded. This gives a limit switched ODE for the sequence of interest as $\mu\to 0$. Here, we study asymptotic stability and establish that the limit points of the switched ODE and the stochastic approximation algorithm coincide as $t\rightarrow \infty$. We thus consider the case where $\mu\to 0$ and $n\to\infty$, however, $\mu n \to\infty$ now. Nevertheless, instead of considering a two-stage limit by first letting $\mu\rightarrow0$ and then $t\rightarrow\infty$, we study $\X^\mu(t+q(\mu))$ and require $q(\mu)\rightarrow\infty$ as $\mu\rightarrow 0$. The following corollary concerns asymptotic stability of the interpolated process.

\vspace{0.1cm}
\begin{Corollary}
\label{cor:weak-convergence}
Let
\begin{equation}
\mathcal{X}^\eta = \lbr
\left[\boldsymbol{x}, r\right]'; \boldsymbol{x} = \mathbf{0}_S, r\in\RR_{[0,\eta)}\rbr.
\end{equation}
Denote by $\left\lbrace q\left(\mu\right)\right\rbrace$ any sequence of real numbers satisfying $q\left(\mu\right)\rightarrow\infty$ as $\mu\rightarrow0$.
Assume $\{\X(n): \mu > 0, n <\infty\}$ is tight or bounded in probability. Then, for each $\eta \geq 0$, there exists $\bar{\gamma}\left(\eta\right)\geq 0$ such that if $\gamma\leq\bar{\gamma}\left(\eta\right)$ in~(\ref{eq:smooth-br}),
\beq{}
\X^\mu(\cdot+q\left(\mu\right))\Rightarrow \mathcal{X}^\eta, \quad \textmd{as $\mu\to 0$}.
\eeq
\begin{IEEEproof}
\newcommand{\wdt}{\hat}
We only give an outline of the proof, which essentially follows from Theorems~\ref{theorem:discrete-continuous} and~\ref{theorem:Stability}. Define $\wdt \X^{\mu}\cd= \X^{\mu}(\cdot + q(\mu))$. Then, it can be shown that $\wdt \X^{\mu}\cd$ is tight. For any $T_1 < \infty$, take a weakly convergent subsequence of $\big\{\wdt \X^{\mu}\cd, \wdt\X^{\mu}(\cdot - T_1)\big\}$. Denote the limit by $\big( \wdt\X\cd,  \wdt\X_{T_1}\cd\big)$. Note that $\wdt\X (0) = \wdt\X_{T_1} (T_1)$. The value of $\wdt\X_{T_1}(0)$ may be unknown, but the set of all possible values of $\wdt\X_{T_1}(0)$ (over all $T_1$  and convergent subsequences) belongs to a tight set. Using this and Theorems~\ref{theorem:discrete-continuous} and~\ref{theorem:Stability}, for any $\varrho > 0$, there exists a $T_\varrho < \infty$ such that for all $T_1 >T_\varrho$, $ d\big(\wdt\X_{T_1}(T_1), {\mathcal X}^\eta\big) \ge 1-\varrho$. This implies that $d\big(\wdt\X (0), {\mathcal X}^\eta\big) \ge 1-\varrho$, and the desired result follows.
\end{IEEEproof}
\end{Corollary}

\subsection{Performance Analysis via Limit Set Characterization}
\label{sec:final_convergence}

The final stage of the proof deals with the analysis of efficiency and tracking properties of the adaptive discrete stochastic optimization algorithm through characterizing the limit set of the switched ODE. The result concerning the tracking capability in Theorem~\ref{theorem:Main} follows directly from Corollary~\ref{cor:weak-convergence}. In what follows, we use this result to conclude efficiency of Algorithm~1 by showing that the empirical sampling distribution $\z(n)$ tracks the set $\Upsilon^\eta(\theta(n))$ (see~(\ref{eq:simplex-epsilon})).

Define the interpolated sequence of iterates $\bar{f}^\mu(t) = \bar{f}(n)$ for $t\in[n\mu,(n+1)\mu)$, and recall the interpolated processes~(\ref{eq:interpol-z-r}). Suppose $\theta(\tau) = \bar{\theta}$ for $\tau\geq 0$. Then, in view of~(\ref{eq:average_1}) and~(\ref{eq:empirical_frequency_local}),
\begin{align}
\label{eq:lim-regret}
r^{\mu}(t) =
\bar{f}^\mu(t) - F_{\min}(\bar{\theta})
 =
 \sum_{i\in\M} z_i^\mu(t)\lb f(i,\bar{\theta}) - F_{\min}(\bar{\theta})\rb,
\end{align}
since $\sum_{i\in\M} z_i^\mu(t) = 1$. On any convergent subsequence $\lbrace\z(n')\rbrace_{n'\geq 0}\rightarrow \boldsymbol{\pi}(\bar{\theta})$, with slight abuse of notation, let $\z^{\mu}(t) = \z(n')$ and $r^{\mu}(t) = r(n')$ for $t\in[n'\mu,n'\mu+\mu)$. This, together with (\ref{eq:lim-regret}), yields
\beq{eq:Lemma_switching_2}
r^{\mu}(\cdot+q(\mu)) \to \sum_{i\in\M}\pi_i(\bar{\theta})\left[f(i,\bar{\theta}) - F_{\min}(\bar{\theta})\right], \;\;\hbox{as }\mu\to 0,
\eeq
since $q(\mu)\to 0$ as $\mu\to 0$. Finally, comparing (\ref{eq:Lemma_switching_2}) with (\ref{eq:simplex-epsilon}) concludes that, for each $\bar{\theta}\in\mathcal{Q}$, $\z^{\mu}(\cdot+q(\mu))$ converges to the $\Upsilon^\eta(\bar{\theta})$ if and only if $r^{\mu}(\cdot+q(\mu)) \leq \eta$ as $\mu\to 0$.
%
%
Combining
this with Corollary~\ref{cor:weak-convergence} completes the proof of the efficiency result in Theorem~\ref{theorem:Main}.

\section{Numerical Examples}
\label{sec:numerical-example}
This section illustrates the performance of Algorithm~1 using the examples in~\cite{andradottir1996global,andradottir1999accelerating}. We start with a static discrete stochastic optimization example, in order to compare Algorithm~1 with two existing algorithms in the literature. We then proceed to the regime-switching framework to illustrate the tracking capability of Algorithm~1.

\subsection{Example 1: Static Discrete Stochastic Optimization}

Consider the following example described in~\cite[Section~4]{andradottir1996global}. Suppose that the demand $Y$ for a particular product has a Poisson distribution with
parameter $\lambda$, i.e., the probability function is give by
\begin{equation}
d\sim f(s;\lambda) = \frac{{\lambda^s \exp(-\lambda)}}{s!}.\nonumber
\end{equation}
The objective is then to find the order size that maximizes the demand probability, subject to the constraint that at most $S$ units can be ordered. This problem can be formulated as a discrete deterministic optimization problem:
\begin{equation}
\label{eq:example-1}
\argmax_{s \in \lbr 0,1,\ldots,S\rbr} \lb f(s;\lambda) = \frac{{\lambda^s \exp(-\lambda)}}{s!}\rb,
\end{equation}
which can be solved analytically. Here, we aim to solve the following stochastic variant: Compute
\begin{equation}
\label{eq:example-2}
\argmin_{s \in \lbr 0,1,\ldots,S\rbr} - \ee\lbr I_{\lbr d = s\rbr}\rbr,
\end{equation}
where $I_{\lbr\cdot\rbr}$ denotes the indicator function, and $d$ is a Poisson distributed random variable with rate $\lambda$. Clearly, problems~(\ref{eq:example-1}) and~(\ref{eq:example-2}) both lead to the same set of global optimizers. This enables us to check the results obtained using Algorithm~1.

We consider the following two cases of the rate parameter $\lambda$ in~(\ref{eq:example-1}): i) $\lambda = 1$, which implies that the set of global optimizers is $\go = \lbr 0,1\rbr$, and ii) $\lambda = 10$, in which case the set of global optimizers is $\go = \lbr 9,10\rbr$. For each case, we further study the effect of the search space size on the performance of Algorithm~1 by considering two instances: i) $S=10$, and ii) $S = 100$. Finally, we compare Algorithm~1 (referred to as AS) with the following two algorithms that have been proposed in the literature:
\begin{enumerate}
    \item[i)] Random search (RS)~\cite{andradottir1996global}: Each iteration of the RS algorithm requires one random number selection, $\mathcal{O}(S)$ arithmetic operations, one comparison and two independent evaluations of the objective function $f_n(s)$.
    \item[ii)] Upper confidence bound (UCB)~\cite{auer2002finite}: Each iteration of the UCB algorithm requires
        $\mathcal{O}(S)$ arithmetic operations, one maximizations and one evaluation of the objective function $f_n(s)$.
\end{enumerate}
Note in comparison that, using $\rho(x)$ as in Remark~\ref{remark-3-1}, the AS algorithm proposed in this paper requires $\mathcal{O}(S)$ arithmetic operations, one random number selection and one evaluation of the objective function $f_n(s)$ at each iteration. Since the problem is static in the sense that $\go$ is fixed for each case, we apply the modifications discussed in Sec.~\ref{sec:static-DSO} to Algorithm~1 and set $\alpha = 0.2$ and $\gamma = 0.01$ in this example.

\begin{table}[!t]
\label{table-both}
  \centering
  \caption{Example 1: Percentage of Cases Where Algorithms Converged to Global Optima $\go(\theta(n))$ in $n$ Iterations}
  $
\begin{array}{l}
\textmd{RS: Random Search Algorithm of~\cite{andradottir1996global}}\\
\textmd{AS: Proposed Adaptive Search in Algorithm~1}\\
\textmd{UCB: Upper Confidence Bound Algorithm of~\cite{auer2002finite}}
\end{array}
$
  \subtable[$\lambda = 1$]{
  \label{table:1}
  \centering
\renewcommand{\arraystretch}{1.2}
\begin{tabular}{cccccccc}
\firsthline
Iteration & \multicolumn{3}{c}{$S = 10$} & & \multicolumn{3}{c}{$S = 100$} \\
\cline{2-4} \cline{6-8}
$n$       &  AS  &      RS     &     UCB     & &  AS   &    RS       &     UCB \\
\hline
10        &  55  &      39     &     86      & &  11   &     6       &      43     \\
50        &  98  &      72     &     90      & &  30   &     18      &      79     \\
100       &  100 &      82     &     95      & &  48   &     29      &      83     \\
500       &  100 &      96     &     100     & &  79   &     66      &      89     \\
1000      &  100 &     100     &     100     & &  93   &     80      &      91     \\
5000      &  100 &     100     &     100     & &  100  &     96      &      99     \\
10000     &  100 &     100     &     100     & &  100  &     100     &      100    \\
\lasthline
\end{tabular}
  }
  \subtable[$\lambda = 10$]{
  \label{table:2}
    \centering
\renewcommand{\arraystretch}{1.2}
\begin{tabular}{cccccccc}
\firsthline
Iteration & \multicolumn{3}{c}{$S = 10$} & & \multicolumn{3}{c}{$S = 100$} \\
\cline{2-4} \cline{6-8}
$n$       & AS & RS & UCB & &  AS   & RS  & UCB \\
\hline
10        &  29  &     14     &    15     & &   7    &     3       &     2      \\
100       &  45  &     30     &    41     & &   16   &     9       &     13     \\
500       &  54  &     43     &    58     & &   28   &     21      &     25     \\
1000      &  69  &     59     &    74     & &   34   &     26      &     30     \\
5000      &  86  &     75     &    86     & &   60   &     44      &     44     \\
10000     &  94  &     84     &    94     & &   68   &     49      &     59     \\
20000     &  100 &     88     &    100    & &   81   &     61      &     74     \\
50000     &  100 &     95     &    100    & &   90   &     65      &     81     \\
\lasthline
\end{tabular}
  }
\end{table}

To give a fair comparison of the three algorithms, we use the iteration number to denote the number of performed simulations. All three algorithms are initialized at state $s(0)$, that is chosen uniformly from $\M$, and move towards $\go$.
Close scrutiny of the results presented in Table~I leads to the following observations: In all three algorithms, the speed of convergence decreases when either $S$ or $\lambda$ (or both) increases. However, the effect of increasing $\lambda$ is more substantial since the objective function values of the worst and best states become closer when $\lambda = 10$.
At a fixed iteration number, higher percentage of cases where a particular method has converged to the global optima indicates convergence at a faster rate. As the results of Table~I show, Algorithm~1 ensures faster convergence to the global optima $\go$ in each case.

\begin{figure}[!t]
\setlength{\abovecaptionskip}{0em}
\setlength{\belowcaptionskip}{0em}
     \centering
     \includegraphics[width=3.6in]{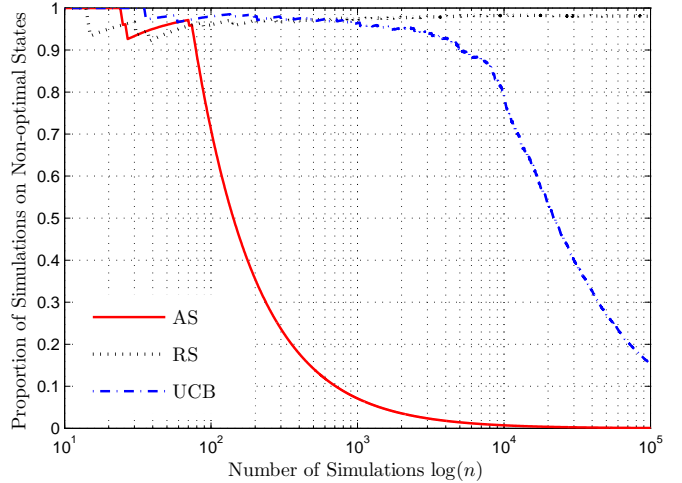}
     \caption{Example 1: Proportion of simulation effort expended on states outside the global optima set ($\lambda = 1$, $S = 100$).}
     \label{fig:1}
\end{figure}

To illustrate superior efficiency of Algorithm~1, we plot the sample path of the time spent simulating states outside the global optima, i.e.,
\begin{equation}
\label{eq:eff}
\textstyle 1- \sum_{i\in\go} z_i(n),
\end{equation}
in Fig.~\ref{fig:1}. This figure corresponds to the case where $\lambda = 1$ and $S=100$ in~(\ref{eq:example-2}). As can be seen, since the RS method randomizes among all states (except the previously sampled state) at each iteration, it spends roughly 98\% of its simulation effort on non-optimal states. Further, the UCB algorithm switches to its exploitation phase after a longer period of exploration as compared to Algorithm~1. Fig~\ref{fig:1} thus indicates that Algorithm~1 guarantees a superior balance between exploration of the search space and exploitation of the collected data as compared to other schemes.

\subsection{Example 2: Regime-Switching Discrete Stochastic Optimization}

Consider the discrete stochastic optimization problem described in Example~1 with the exception that now $\lambda(\theta(n))$ jump changes between 1 and 10 according to a slow Markov chain $\lbr\theta(n)\rbr$ with state space $\mathcal{Q} = \lbr 1,2\rbr$, and transition probability matrix
\begin{equation}
\label{eq:Markov-example}
P^\vare = I+\vare Q,\quad Q = \left[\begin{matrix} -0.5 & 0.5 \\ 0.5 & -0.5 \end{matrix}\right].
\end{equation}
More precisely, $\lambda(1) = 1$ and $\lambda(2) = 10$. Assuming $S = 10$, $\go(1) = \lbr 0,1\rbr$ and $\go(2) = \lbr 9,10\rbr$. Then, the discrete stochastic optimization problem is given by~(\ref{eq:example-2}), where
\begin{equation}
d\sim f(s,
i;\lambda) = \frac{{\lambda^s(
i) \exp(-\lambda(i
))}}{s!}, \ i=1,2.
\end{equation}
In the rest of this section, we assume $\gamma = 0.1$, and $\mu = \vare = 0.01$. Further, we shall use an adaptive variant of RS, studied in~\cite{YKI04}, and an adaptive variant of UCB both with constant step-sizes $\mu = 0.01$ to compare with algorithm~1.

\begin{figure}[!t]
\setlength{\abovecaptionskip}{0em}
\setlength{\belowcaptionskip}{0em}
     \centering
     \includegraphics[width=3.6in]{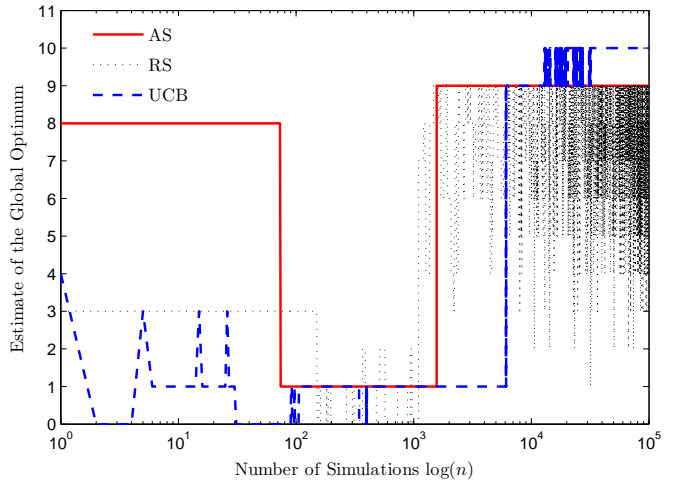}
     \caption{Example 2: Sample path of the estimate of the global optima when $S = 10$ and the global optima evolve with time. For $0\leq n < 10^3$, the global optima set is $\lbr 0,1\rbr$. For  $10^3\leq n<10^5$, the global optima set is $\lbr9,10\rbr$.}
     \label{fig:2}
\end{figure}

\begin{figure}[!t]
\setlength{\abovecaptionskip}{0em}
\setlength{\belowcaptionskip}{0em}
     \centering
     \includegraphics[width=3.6in]{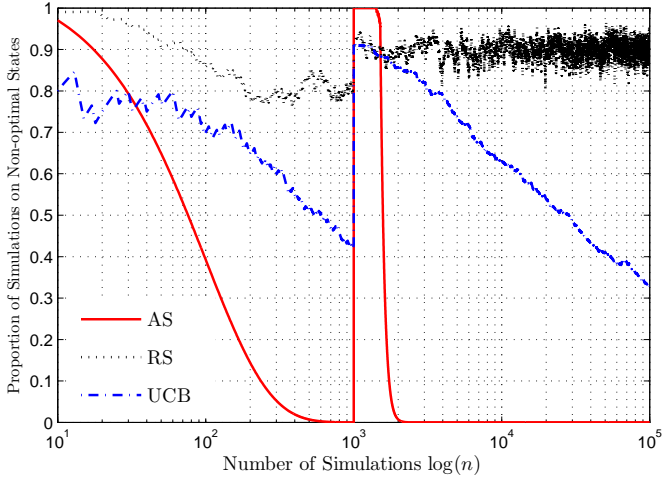}
     \caption{Example 2: Proportion of simulation effort expended on non-optimal states when $S = 10$ and the global optima evolve with time. For $0\leq n < 10^3$, the global optima set is $\lbr 0,1\rbr$. For  $10^3\leq n<10^5$, the global optima set is $\lbr9,10\rbr$.}
     \label{fig:3}
\end{figure}

Fig.~\ref{fig:2} shows tracking capability of Algorithm~1 when the Markov chain $\lbr \theta(n)\rbr$ undergoes a jump from $\theta = 1$ to $\theta = 2$ at $n = 10^3$. As can be seen, contrary to the RS algorithm, both AS and UCB methods properly track the changes; however, AS is more agile. Superior performance of the AS algorithm is further verified in Fig.~\ref{fig:3} which shows how the simulation effort on non-optimal states evolves as the rate parameter $\lambda$
jump changes. Fig.~\ref{fig:3} thus confirms that the superior balance between exploration and exploitation properly responds to the regime switching.

\begin{figure}[!t]
\setlength{\abovecaptionskip}{0em}
\setlength{\belowcaptionskip}{0em}
     \centering
     \includegraphics[width=3.6in]{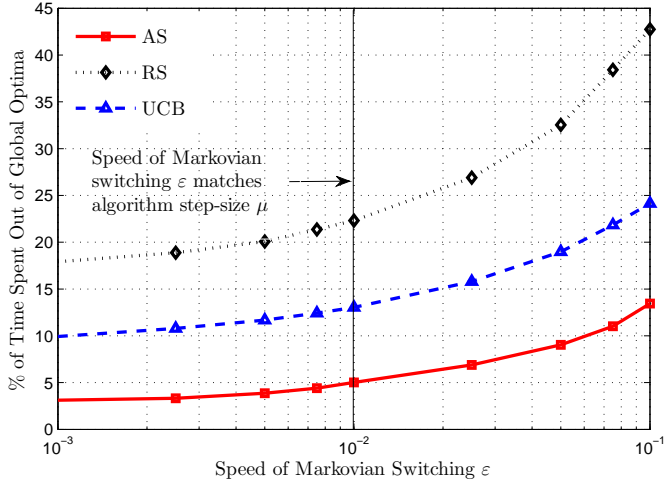}
     \caption{Example 2: Proportion of time the estimate of global optima spends out of the global optima set $\go(\theta(n))$ versus the speed of Markovian switching of the set of global optimizers~($S = 10$).}
     \label{fig:4}
\end{figure}

Fig.~\ref{fig:4} illustrates the efficiency~(\ref{eq:eff}) of the AS algorithm for several values of $\vare$. Note that $\vare$ represents the speed of Markovian switching. Each point on the graph is an average over 100 independent runs of $10^6$ iterations of the algorithms when~(\ref{eq:Markov-example}) is adopted as the transition matrix of $\lbr\theta(n)\rbr$. As expected, the percentage of samples taken from the set of global optima increases for all methods as the speed of time variations decreases; however, superior efficiency of the AS algorithm is clearly evident.

\def\eem{\ee^\mu_t}

\def\W{\boldsymbol{W}}
\begin{figure*}[!hb]
\hrulefill
\vspace*{-2pt}
\setcounter{MYtempeqncnt}{\value{equation}}
\setcounter{equation}{41}
\begin{align}
\label{eq:A1-tight-first}
\eem \left\| \mu \sum^{(t+u)/\mu-1}_{k=t/\mu} \boldsymbol{A}_k\left(s(k)\right) - \X(k)\right\|^2
&= \mu^2 \sum^{(t+u)/\mu-1}_{\tau=t/\mu} \sum^{(t+u)/\mu-1}_{\kappa=t/\mu} \eem \lbr \lb\boldsymbol{A}_\tau\left(s(\tau)\right) - \X(\tau)\rb' \lb\boldsymbol{A}_\tau\left(s(\kappa)\right) - \X(\kappa)\rb\rbr \nonumber\\
&\leq K\mu^2 \left(\frac{t+u}{\mu} - \frac{t}{\mu}\right)^2 = O\left(u^2\right)
\end{align}
%
\setcounter{equation}{\value{MYtempeqncnt}}
\vspace*{-17pt}
\end{figure*}

\section{Conclusion}
\label{sec:conclusion}

This paper has considered regime-switching discrete stochastic optimization problems where the underlying time variations, e.g., in the profile of the stochastic behavior of the system or the objective function, can be captured by the sample path of a slow discrete time Markov chain. We proposed a class of adaptive search algorithms that prescribes how to iteratively sample states from the search space. The proposed scheme is a constant step-size stochastic approximation algorithm that updates beliefs about the objective function values, accompanied by an adaptive sampling strategy of best-response type. The convergence analysis proved that, if the underlying time variations occur on the same timescale as
the stochastic approximation algorithm, the algorithm will properly track the randomly switching set of global minima.
Further, the proposed scheme ensures ``most'' of the simulation effort is spent on the global minima. It thus can be deployed as an on-line control mechanism to enable self-configuration of large scale stochastic systems. The main features of the proposed adaptive discrete stochastic optimization algorithm include: 1) it allows time correlation in the sampled data; 2) it tracks time varying optima when the parameters underlying the stochastic optimization problem evolve over time; 3) in contrast to the case where the time variations occur on a slower timescale as the adaptive search algorithm and trackability is trivial, it tracks time variations of the global optima even when such variations occur on the same timescale as the updates of the proposed algorithm.
Numerical examples
illustrated the
trade-off between efficiency (the number of executed simulations) and the convergence speed, as compared with the existing random search and pure exploration methods.


\appendices

\section{Proof of Theorem~\ref{theorem:discrete-continuous}}
\label{sec:proof_1}

We first prove tightness of the interpolated process $\X^\mu\cd$. Consider the sequence $\lbr \X(n)\rbr$, defined in~(\ref{eq:X_def}). In view of the boundedness of the objective function, and by virtue of H\"older's and Gronwall's inequalities, for any $0<T_1<\infty$,
\beq{}
\sup_{k \le T_1/\mu} \ee
  \left\| \X(k)\right\|^2<\infty,
\eeq
where in the above and hereafter $\|\cdot\|$ denotes the Euclidean norm and $t/\mu$ is understood to be the integer part of $t/\mu$ for each $t>0$. Next, considering the interpolated process $\X^\mu\cd$ (defined in (\ref{eq:interpolations})) and the recursion~(\ref{eq:X-recursion}), for any $t,u>0$, $\delta>0$, and $u<\delta$, it can be verified that
\beq{}
\label{eq:A1-recursion}
\begin{split}
&\X^\mu(t+u) - \X^\mu(t) = \mu \sum^{(t+u)/\mu-1}_{k=t/\mu} \lb\boldsymbol{A}_k\left(s(k)\right) - \X(k)\rb\\
& \hspace{1.2cm} + \sum^{(t+u)/\mu-1}_{k=t/\mu} \left[
\begin{array}{c}
\boldsymbol{F}(\theta(k)) - \boldsymbol{F}(\theta(k+1))\\
F_{\min}(\theta(k)) - F_{\min}(\theta(k+1))
\end{array}
\right],
\end{split}
\eeq
where $\boldsymbol{A}_k\left(s(k)\right)$ is defined in~(\ref{eq:A1-def}).
%
%
%
%
Consequently, using the parallelogram law,
\begin{equation}
\label{eq:A1-2}
\begin{split}
&\eem \left\|\X^\mu(t+u) - \X^\mu(t)\right\|^2 \\
&\;\;\leq 2\eem \left\| \mu \sum^{(t+u)/\mu-1}_{k=t/\mu} \boldsymbol{A}_k\left(s(k)\right) - \X(k)\right\|^2\\
&\;\;+ 2\eem \left\|\sum^{(t+u)/\mu-1}_{k=t/\mu}  \left[
\begin{array}{c}
\boldsymbol{F}(\theta(k)) - \boldsymbol{F}(\theta(k+1))\\
F_{\min}(\theta(k)) - F_{\min}(\theta(k+1))
\end{array}
\right] \right\|^2,
\end{split}
\end{equation}
where $\mathbb{E}^\mu_t$ denotes the $\sigma$-algebra generated by the $\mu$-dependent past data up to time $t$. By virtue of the tightness criteria~\cite[Theorem 3, p. 47]{kushner1984approximation} or~\cite[Chapter 7]{KY03}, it suffices to verify
\beq{eq:A1-3}
\disp\lim_{\delta\to 0}
\limsup_{\mu\to 0}
\lbr
\mathbb{E}
\lb\sup_{0\le u \le \delta} \mathbb{E}^\mu_t\left \|
\X^\mu(t+u) - \X^\mu(t)
\right\|^2
\rb
\rbr  = 0.
\eeq
As for the first term on the r.h.s. of~(\ref{eq:A1-2}), noting the boundedness of objective function, we obtain: see~(\ref{eq:A1-tight-first}) at the bottom of the page. We then concentrate on the second term on the r.h.s. of~(\ref{eq:A1-2}).
Note that, for sufficiently small positive $\mu$, if $Q$ is irreducible, then so is $I + \mu Q$.
Thus, for sufficiently large~$k$,
$\| (I + \mu Q)^k - \one \nu_\mu\|_M \le \lambda_c^k$
for some $0<\lambda_c<1$, where $\nu_\mu$ denotes the row vector of
stationary distribution associated with the transition matrix $I+ \mu
Q$, $\one$ denotes the column vector of ones, and $\|\cdot\|_M$ represents any matrix norm.
The essential feature involved in the second term in~(\ref{eq:A1-2})
is the difference of the transition probability matrix of the form
$(I+ \mu Q) ^{k-(t/\mu)}- (I+\mu Q)^{k+1 - (t/\mu)}$. However, it can be
seen that
\bea \ad\!\!\!\! (I+ \mu Q)^{k -(t/\mu)} -(I+\mu Q)^{k+1-(t/\mu)}
 \\
\ad  = -\mu Q[(I+\mu Q)^{k-(t/\mu)} - \one \nu_\mu ] .\eea
%
This in turn implies that
\setcounter{equation}{42}
\begin{equation}
\label{eq:A1-tight-second}
\begin{split}
\eem &\left\|\sum^{(t+u)/\mu-1}_{k=t/\mu} \left[
\begin{array}{c}
\boldsymbol{F}(\theta(k)) - \boldsymbol{F}(\theta(k+1))\\
F_{\min}(\theta(k)) - F_{\min}(\theta(k+1))
\end{array}
\right] \right\|^2\\
\aad\leq K O(\mu) \sum^{(t+u)/\mu-1}_{k=t/\mu} \lambda_c^k \le O(\mu) \sum^\infty_{k=1} \lambda^k_c  = O(\mu).
\end{split}
\end{equation}
Finally, combining~(\ref{eq:A1-tight-first}) and~(\ref{eq:A1-tight-second}), the tightness criteria~(\ref{eq:A1-3}) is verified. Therefore, $\X^\mu\cd$ is tight in $D([0,\infty]:\RR^{S+1})$. In view of~\cite[Proposition 4.4]{YKI04}, $\theta^\mu\cd$ is also tight and $\theta^\mu\cd\Rightarrow \theta\cd$ such that $\theta\cd$ is a continuous time Markov chain with generator $Q$; see~(A1). As the result, the pair $(\X^\mu\cd,\theta^\mu\cd)$ is tight in $D\big([0,\infty]:\RR^{S+1}\times\mathcal{Q}\big)$.

Using Prohorov's theorem~\cite{KY03}, one can extract a convergent subsequence. For notational simplicity, we still denote the subsequence by $\X^\mu\cd$ with limit $\X\cd$. By the Skorohod representation theorem~\cite{KY03}, and with a slight abuse of notation, $\X^\mu\cd \rightarrow \X\cd$ in the sense of w.p.1 and the convergence is uniform on any compact interval. We now proceed to characterize the limit $\X\cd$ using martingale averaging methods.

First, we demonstrate that the last term in~\eqref{eq:A1-recursion}
contributes nothing to the limit differential equation.
We aim to show
\begin{equation}\label{ef-fmin}
\hspace{-0.3cm}\begin{array}{ll} \ad
\lim_{\mu\to 0}\ee h\left(\X^\mu(t_\iota),\theta^\mu(t_\iota):\iota\leq\kappa_0\right)
\\
\aad \times E^\mu_t \Bigg[ \sum^{ (t+u)/\mu -1}_{k=t/\mu} \left[
\begin{array}{c}
\boldsymbol{F}(\theta(k)) - \boldsymbol{F}(\theta(k+1))\\
F_{\min}(\theta(k)) - F_{\min}(\theta(k+1))
\end{array}
\right]=0. \end{array}\nonumber
\end{equation}
This directly follows from an argument similar to the one used in~(\ref{eq:A1-tight-second}).

To obtain the desired limit, it will be proved that the limit $(\X\cd, \theta\cd)$ is the solution of the martingale problem with operator $\mathcal{L}$ defined as follows: For all $i\in\mathcal{Q}$,
\beq{eq:L1-operator}
\begin{split}
&\mathcal{L} y(x,i) = \nabla_x y'(x,i)\left[\boldsymbol{G}(x,i)-x\right] + Q y(x,\cdot)(i),\\
&\hspace{1.3cm}Q y(x,\cdot)(i) = \sum_{j\in \mathcal{Q}} q_{i j} y(x,j),
\end{split}
\eeq

\begin{figure*}[!hb]
\setcounter{MYtempeqncnt}{\value{equation}}
\setcounter{equation}{46}
\hrulefill
\vspace*{-2pt}
\begin{align}
\label{eq:A-4}
&\lim_{\mu\to 0}
\ee h(\X^\mu(t_\iota), \theta^\mu(t_\iota): \iota\le \kappa_0)\Bigg[\sum^{ t+u}_{\ell:\ell \delta_\mu = t}\lb y\left(\X(\ell n_\mu), \theta(\ell n_\mu + \mu)\right) - y\left(\X(\ell n_\mu),\theta(\ell n_\mu)\right)\rb\Bigg]\nonumber\\
 &\;\;=  \lim_{\mu\to 0}
\ee h\big(\X^\mu(t_\iota), \theta^\mu(t_\iota): \iota\le \kappa_0\big)\nonumber\\
&\hspace{1.5cm}\times\Bigg[\sum^{ t+u}_{\ell:\ell \delta_\mu =t }
\sum^\Theta_{i_0=1} \sum^\Theta_{j_0=1}
\sum^{\ell n_\mu+n_\mu -1}_{k=\ell n_\mu}\lb y\left(\X(\ell n_\mu),j_0\right) \mathds{P}( \theta(k+1)=j_0| \theta(k) =i_0) - y(\X(\ell n_\mu),i_0)\rb I_{\{\theta(k) =i_0\}}\Bigg]\nonumber\\
&\;\; =\ee
h\big(\X(t_\iota), \theta(t_\iota): \iota\le \kappa_0\big)
\left[ \int^{t+u}_t
Q y(\X(v), \theta(v)) dv\right]
\end{align}
\setcounter{equation}{\value{MYtempeqncnt}}
\vspace*{-17pt}
\end{figure*}

\vspace{-0.1cm}
\noindent
and, for each $i\in\mathcal{Q}$, $y(\cdot, i): \RRk \mapsto \rr$ with $y(\cdot, i) \in C^1_0$ ($C^1$ function with compact support). Further, $\nabla_x y(x,i)$ denotes the gradient of $y(x,i)$ with respect to $x$, and $\boldsymbol{G}(\cdot,\cdot)$ is defined in (\ref{eq:G}). Using an argument similar to~\cite[Lemma 7.18]{yin1998continuous}, one can show that the martingale problem associated with the operator $\mathcal{L}$ has a unique solution. Thus, it remains to prove that the limit $(\X\cd, \theta\cd)$ is the solution of the martingale problem. To this end, it suffices to show that, for any positive arbitrary integer $\kappa_0$, and for any $t,u>0$, $0<t_\iota\leq t$ for all $\iota\leq\kappa_0$, and any bounded continuous function $h(\cdot,i)$ for all $i\in\mathcal{Q}$,
\beq{eq:martingale}
\begin{split}
&\ee h(\X(t_\iota), \theta(t_\iota): \iota\le \kappa_0)\\
&\quad\times\bigg[y\left(\X(t+u),\theta(t+u)\right) - y\left(\X(t),\theta(t)\right) \\
&\qquad\;\;- \int_t^{t+u} \mathcal{L} y\left(\X(v),\theta(v)dv\right)\bigg] = 0.
\end{split}
\eeq
To verify (\ref{eq:martingale}), we work with $(\X^\mu\cd,\theta^\mu\cd)$ and prove that the above equation holds as $\mu\to 0$.

\begin{figure*}[!hb]
\vspace*{-2pt}
\setcounter{MYtempeqncnt}{\value{equation}}
\setcounter{equation}{47}
\begin{align}
\label{eq:A-5}
&\lim_{\mu\to 0}
 \ee h\left(\X^\mu(t_\iota), \theta^\mu( t_\iota):\iota\leq\kappa_0\right)\Bigg[\sum^{ t+u}_{\ell:\ell \delta_\mu =t }\left[y\left(\X(\ell n_\mu + n_\mu), \theta(\ell n_\mu + n_\mu)\right) - y(\X\left(\ell n_\mu), \theta(\ell n_\mu + n_\mu)\right)\right]\Bigg]\nonumber\\
&\;\; =\lim_{\mu\to 0} \ee h\left(\X^\mu(t_\iota), \theta^\mu( t_\iota):\iota\leq\kappa_0\right)
 \Bigg[\sum^{ t+u}_{\ell:\ell \delta_\mu =t }\left[y\left(\X(\ell n_\mu + n_\mu), \theta(\ell n_\mu)\right) - y\left(\X(\ell n_\mu), \theta(\ell n_\mu)\right)\right]\Bigg]\nonumber\\
&\;\; = \lim_{\mu\to 0}\ee h\left(\X^\mu(t_\iota),\theta^\mu( t_\iota):\iota\leq\kappa_0\right)
 \Bigg[\sum^{ t+u}_{\ell:\ell \delta_\mu =t }
 \nabla^\prime_{_{\fhat}} y\left(\X(\ell n_\mu), \theta(\ell n_\mu)\right)\left[\fhat(\ell n_\mu + n_\mu)- \fhat(\ell n_\mu)\right] \nonumber\\
 &\hspace{6.2cm}+ \nabla^\prime_{_{r}} y\left(\X(\ell n_\mu), \theta(\ell n_\mu)\right)\left[r(\ell n_\mu + n_\mu)- r(\ell n_\mu)\right]\Bigg]
 \nonumber\\
&\;\;= \lim_{\mu\to 0}\ee h\left(\X^\mu(t_\iota),\theta^\mu( t_\iota):\iota\leq\kappa_0\right)\nonumber\\
 &\hspace{1.5cm}\times\Bigg[\sum^{ t+u}_{\ell:\ell \delta_\mu =t }
 \delta_\mu \nabla^\prime_{_{\fhat}} y\left(\X(\ell n_\mu), \theta(\ell n_\mu)\right)\Bigg[\frac{1}{n_\mu}
 \sum^{\ell n_\mu+n_\mu-1}_{k=\ell n_\mu}
 \lb \hat{\boldsymbol{g}}_{_{\theta(\ell n_\mu)}}\big(s(k),\fhat(k)\big) - \bf(\theta(\ell n_\mu))\rb
 - \frac{1}{n_\mu} \sum^{\ell n_\mu + n_\mu-1}_{k=\ell n_\mu}\fhat(k) \Bigg]\nonumber\\
 &\hspace{3cm} + \delta_\mu \nabla^\prime_{_{r}} y\left(\X(\ell n_\mu), \theta(\ell n_\mu)\right)\Bigg[\frac{1}{n_\mu}
 \sum^{\ell n_\mu+n_\mu-1}_{k=\ell n_\mu} \lb
 f_k\left(\boldsymbol{s}(k)\right) - F_{\min}(\theta(\ell n_\mu))\rb
 - \frac{1}{n_\mu} \sum^{\ell n_\mu + n_\mu-1}_{k=\ell n_\mu} r(k) \Bigg]\Bigg]
\end{align}
\setcounter{equation}{\value{MYtempeqncnt}}
\vspace*{-17pt}
\end{figure*}

By the weak convergence of $(\X^\mu\cd,\theta^\mu\cd)$ to $(\X\cd,\theta\cd)$ and Skorohod representation, it can be
seen that
\beq{}
\begin{split}
&\ee
h\left(\X^\mu(t_\iota), \theta^\mu(t_\iota): \iota\le \kappa_0\right)\\
&\quad\times\left[ \left(\X^\mu(t+u),\theta^\mu(t+u)\right)- \left(\X^\mu(t),\theta^\mu(t)\right)\right] \\
&\to \ee h\left(\X(t_\iota), \theta(t_\iota): \iota\le \kappa_0\right)\\
&\qquad\times\left[ \left(\X(t+u),\theta(t+u)\right)- \left(\X(t),\theta(t)\right)\right].
\end{split}
\nonumber
\eeq
Now, choose a sequence of integers $\lbrace n_\mu\rbrace$ such that $n_\mu\rightarrow \infty$ as $\mu\rightarrow 0$, but $\delta_\mu = \mu n_\mu \rightarrow 0$, and Partition $[t,t+u]$ into subintervals of length $\delta_\mu$.
Then,
\begin{align}
\label{eq-split}
y\big(&\X^\mu(t+u), \theta^\mu(t+u)\big) - y\left(\X^\mu(t), \theta^\mu(t)\right)\nonumber\\
 &= \sum^{ t+u}_{\ell: \ell \delta_\mu =t }
\Big[ y\left(\X(\ell n_\mu + n_\mu), \theta(\ell n_\mu + n_\mu)\right)\nonumber\\
&\hspace{1.5cm}- y\left(\X(\ell n_\mu), \theta(\ell n_\mu)\right)\Big]\nonumber\\
&= \sum^{ t+u}_{\ell: \ell \delta_\mu =t }\Big[y\left(\X(\ell n_\mu + n_\mu), \theta(\ell n_\mu + n_\mu)\right) \nonumber\\
&\hspace{1.5cm}-y(\X(\ell n_\mu), \theta(\ell n_\mu + n_\mu))\Big]\nonumber\\
&+ \sum^{ t+u}_{\ell:\ell \delta_\mu =t }\Big[y(\X(\ell n_\mu), \theta(\ell n_\mu + n_\mu)) \nonumber\\
&\hspace{1.5cm}-y(\X(\ell n_\mu), \theta(\ell n_\mu))\Big],
\end{align}
where $\sum^{ t+u}_{\ell:\ell \delta_\mu =t }$ denotes the sum over $\ell$ in the range $t\leq \ell\delta_\mu\leq t+u$.

\def\nfy{\nabla'_{_{\fhat}} y}
\def\nfz{\nabla'_{_{\z}} y}
\def\nfm{\nabla'_{_{r}} y}
\def\nfyy{\nabla_{_{\fhat}} y}
\def\nfzz{\nabla_{_{\z}} y}
\def\nfmm{\nabla_{_{r}} y}
First, we consider the second term on the r.h.s. of~\eqref{eq-split}: see~(\ref{eq:A-4}) at the bottom of the next page. As for the first term on the r.h.s. of~\eqref{eq-split}: see~(\ref{eq:A-5}) at the bottom of the next page, where $\nabla_{_{\boldsymbol{x}}} y$ denotes the gradient column vector with respect to vector $\boldsymbol{x}$, $\nabla'_{_{\boldsymbol{x}}} y$ represents its transpose, and $\hat{\boldsymbol{g}}_\theta(\cdot,\cdot)$ denotes the vector $\hat{\boldsymbol{g}}(\cdot,\cdot)$ in~(\ref{eq:A1-def}) when $\theta(k) = \theta$ is held fixed.
The rest of the proof is divided into two steps, each concerning one of the two terms
in~(\ref{eq:A-5}).
For notational simplicity, we shall write
$\nabla_{_{\fhat}}y(\X(\ell n_\mu), \theta(\ell n_\mu))$, and
$\nabla_{_{r}}y(\X(\ell n_\mu), \theta(\ell n_\mu))$
as $\nfyy$, and $\nfmm$, respectively.

\subsubsection*{Step 1}
We start by looking at
\setcounter{equation}{48}
\begin{align}
&\hspace{-0.3cm}\lim_{\mu\to 0}\ee h\left(\X^\mu(t_\iota),\theta^\mu(t_\iota):\iota\leq\kappa_0\right) \Bigg[ \sum^{ t+u}_{\ell:\ell \delta_\mu =t } \delta_\mu  \nfy \nonumber\\
&\hspace{0.5cm}\times\Bigg[\frac{1}{n_\mu}\sum_{k=\ell n_\mu}^{\ell n_\mu+n_\mu-1} \lb\hat{\boldsymbol{g}}_{_{\theta(\ell n_\mu)}}\big(s(k),\fhat(k)\big) - \bf(\theta(\ell n_\mu))\rb \Bigg]\Bigg]\nonumber\\
&\hspace{-0.3cm} = \lim_{\mu\to 0}\ee h\left(\X^\mu(t_\iota),\theta^\mu(t_\iota):\iota\leq\kappa_0\right) \Bigg[ \sum^{ t+u}_{\ell:\ell \delta_\mu =t } \delta_\mu \nfy \nonumber\\
&\hspace{0.3cm}\times\Bigg[ - \bf(\theta(\ell n_\mu)) +  \frac{1}{n_\mu} \sum_{k=\ell n_\mu}^{\ell n_\mu+n_\mu-1} \sum_{\check \theta=1}^\Theta  \\
&\hspace{0.9cm}\Bigg[ \sum_{\theta = 1}^\Theta \ee_{\ell n_\mu} \hat{\boldsymbol{g}}_{_{\theta(\ell n_\mu)}}\big(s(k),\fhat(k)\big)
\ee_{\ell n_\mu} I_{\lbr \theta(k) = \theta |\theta(\ell n_\mu) = \check \theta\rbr}\Bigg] \Bigg] \Bigg]\nonumber
\end{align}
We
concentrate on the term involving the Markov chain $\theta(k)$. Note that for large $k$ with $\ell n_\mu \le k\le \ell n_\mu +n_\mu$ and $k-\ell n_\mu\to \infty$,
by~\cite[Proposition 4.4]{YKI04}, for some $\hat{k}_0 >0$,
\beq{eq:A1-17}
\begin{split}
(I+\mu Q)^{k-\ell n_\mu} &= Z((k-\ell n_\mu) \mu) \\
&\quad+ O( \mu + \exp(- \hat{k}_0 ( k-\ell n_\mu)),\\
 \frac{dZ(t)}{dt} &=Z(t) Q , \ Z(0)=I.
\end{split}\nonumber
\eeq
For $ \ell n_\mu \le k \le \ell n_\mu + n_\mu$, letting $\mu \ell n_\mu \to u$ yields that $(k -\ell n_\mu) \mu \to 0$ as $\mu\to 0$. For such $k$, $Z((k-\ell n_\mu)\mu) \to I$. Therefore, by the boundedness of $\hat{\boldsymbol{g}}\big(s(k),\fhat(k)\big)$, it follows that, as $\mu\to 0$,
\beq{eq:A1-18}
\begin{split}
&\frac{1}{n_\mu}
 \sum^{\ell n_\mu+n_\mu-1}_{k=\ell n_\mu}
\left\|\ee_{\ell n_\mu} \hat{\boldsymbol{g}}_{_{\theta(\ell n_\mu)}}\big(s(k),\fhat(k)\big)\right\|\\
& \hspace{0.8cm}\times
\Big|\ee_{\ell n_\mu}
\lb I_{\{\theta(k) =\theta\}}| I_{\{ \theta(\ell n_\mu)=\check \theta\}}\rb -I_{\{\theta(\ell n_\mu)=\check \theta\}} \Big|
 \to 0.
\end{split}\nonumber
\eeq
%
%
Therefore,
\begin{align}
\label{eq:A1-12}
&\lim_{\mu\to 0}\ee h\left(\X^\mu(t_\iota),\theta^\mu(t_\iota):\iota\leq\kappa_0\right) \Bigg[ \sum^{ t+u}_{\ell:\ell \delta_\mu =t } \delta_\mu\nfy \nonumber\\
&\hspace{0.4cm}\times\Bigg[\frac{1}{n_\mu}\sum_{k=\ell n_\mu}^{\ell n_\mu+n_\mu-1}
\lb \hat{\boldsymbol{g}}_{_{\theta(\ell n_\mu)}}\big(s(k),\fhat(k)\big) - \bf(\theta(\ell n_\mu))\rb\Bigg]\Bigg]\nonumber\\
&=\lim_{\mu\to 0}\ee h\left(\X^\mu(t_\iota),\theta^\mu(t_\iota):\iota\leq\kappa_0\right)\nonumber\\
 &\hspace{0.4cm}\times\Bigg[ \sum^{ t+u}_{\ell:\ell \delta_\mu =t } \delta_\mu \nfy \sum_{\check \theta=1}^\Theta I_{\lbr \theta(\ell n_\mu) = \check \theta\rbr} \\
&\hspace{1cm}\times\bigg[ - \bf(\check \theta) + \frac{1}{n_\mu}\sum_{k=\ell n_\mu}^{\ell n_\mu+n_\mu-1} \ee_{\ell n_\mu} \hat{\boldsymbol{g}}_{\check \theta}\big(s(k),\fhat(k)\big) \bigg] \Bigg].\nonumber
\end{align}
It is more convenient to work with the individual elements of $\hat{\boldsymbol{g}}_{\check \theta}(\cdot,\cdot)$. Substituting for the $i$-th element from~(\ref{eq:A1-def}) in~(\ref{eq:A1-12}) results
\begin{align}
\label{eq:A1-13}
&\lim_{\mu\to 0}\ee h\left(\X^\mu(t_\iota),\theta^\mu(t_\iota):\iota\leq\kappa_0\right) \nonumber\\
&\hspace{0.5cm}\times\Bigg[ \sum^{ t+u}_{\ell:\ell \delta_\mu =t } \delta_\mu \nfy  \sum_{\check \theta=1}^\Theta I_{\lbr \theta(\ell n_\mu) = \check \theta\rbr} \Bigg[ -F(i, \check \theta)\nonumber\\
& \hspace{1cm}+ \frac{1}{n_\mu} \sum_{k=\ell n_\mu}^{\ell n_\mu+n_\mu-1} \ee_{\ell n_\mu} \Bigg\lbrace \frac{f_k(i)}{{\bbr_i\big(\fhat(k)+\bf(\check \theta)\big)}}\cdot I_{\lbr s(k) = i\rbr}\Bigg\rbrace  \Bigg]\Bigg]\nonumber\\
&=\lim_{\mu\to 0}\ee h\left(\X^\mu(t_\iota),\theta^\mu(t_\iota):\iota\leq\kappa_0\right)\Bigg[ \sum^{ t+u}_{\ell:\ell \delta_\mu =t } \delta_\mu \nfy \\
&\hspace{0.7cm}\times   \sum_{\check \theta=1}^\Theta \Bigg[-F(i,\check \theta) + \frac{1}{n_\mu}  \sum_{k=\ell n_\mu}^{\ell n_\mu+n_\mu-1} \ee_{\ell n_\mu}  f_k(i)\Bigg] I_{\lbr \theta(\ell n_\mu) = \check \theta\rbr} \Bigg].\nonumber
\end{align}
In (\ref{eq:A1-13}), we used $\ee_{\ell n_\mu} I_{\lbr s(k) = i\rbr} = \bbr_i\big(\f(\ell n_\mu)\big)$ since $s(k)$ is chosen according to the smooth best-response strategy $\br(\f(k))$; see Step~1) in Algorithm~1. Note that $f_k(i)$ is still time-dependent due to the presence of noise in the simulation data.
%
%
%
Note further that $\theta(\ell n_\mu) = \theta^\mu( \mu \ell n_\mu)$. In light of (C1)--(C2), by the weak convergence of $\theta^\mu\cd$ to $\theta\cd$, the Skorohod representation, and using $\mu \ell n_\mu\to u$, it can be shown for the second term in~(\ref{eq:A1-13}) that, as $\mu\to 0$,
\beq{eq:A-14}
\begin{split}
& \sum_{\check \theta=1}^\Theta \frac{1}{n_\mu} \sum_{k=\ell n_\mu}^{\ell n_\mu + n_\mu -1} \ee_{\ell n_\mu} f_k(i)  I_{\lbr \theta^\mu(\mu\ell n_\mu) = \check \theta\rbr}  \\
&\quad\to \sum_{\check \theta=1}^\Theta F(i,\check \theta) I_{\{ \theta(u) = \check \theta\}} = F(i,\theta(u))  \quad \textmd{in probability.}\\
\end{split}
\eeq
Using a similar argument for the first term in~(\ref{eq:A1-13}) yields
\begin{align}
\label{eq:A-11}
&\lim_{\mu\to 0}\ee h\left(\X^\mu(t_\iota),\theta^\mu(t_\iota):\iota\leq\kappa_0\right) \Bigg[ \sum^{ t+u}_{\ell:\ell \delta_\mu =t } \delta_\mu\nfy \nonumber\\
&\hspace{0.4cm}\times\Bigg[\frac{1}{n_\mu}\sum_{k=\ell n_\mu}^{\ell n_\mu+n_\mu-1}
\lb \hat{\boldsymbol{g}}_{_{\theta(\ell n_\mu)}}\big(s(k),\fhat(k)\big) - \bf(\theta(\ell n_\mu))\rb\Bigg]\Bigg]\nonumber\\
&\hspace{1cm}\to \mathbf{0}_S \quad \hbox{ as } \mu \to 0.
\end{align}
By using the technique of stochastic approximation (see, e.g., \cite[Chapter 8]{KY03}), it
can be shown that
%
\begin{align}
\label{eq:CIR}
&
\ee h\big(\X^\mu(t_\iota), \theta^\mu(t_\iota): \iota\le \kappa_0\big)\Bigg[\sum_{\ell:\ell \delta_\mu=t}^{t+u} \delta_\mu \nfy \Bigg[\frac{1}{n_\mu} \sum_{k = \ell n_\mu}^{\ell n_\mu+n_\mu-1} \fhat(k)\Bigg]\Bigg]\nonumber\\
& \to \ee h\big(\X(t_\iota), \theta(t_\iota): \iota\le \kappa_0\big)\nonumber\\
&\qquad\times \left[ \int^{t+u}_t \nfy(\X(v),\theta(v)) \fhat(v) dv\right]\;\;\textmd{as}\;\mu\to 0.
\end{align}

\subsubsection*{Step 2}
Next, we concentrate on the second term in~(\ref{eq:A-5}). By virtue of the boundedness of $f_k(s(k))$, and using a similar argument as in Step 1,
%
\begin{align}
\label{eq:A1-19}
&\lim_{\mu\to 0}\ee h\left(\X^\mu(t_\iota),\theta^\mu(t_\iota):\iota\leq\kappa_0\right) \Bigg[ \sum^{ t+u}_{\ell:\ell \delta_\mu =t } \frac{\delta_\mu}{n_\mu} \nfm  \sum_{k=\ell n_\mu}^{\ell n_\mu+n_\mu-1} \nonumber\\
&\hspace{0.5cm} \Bigg[\sum_{\check \theta =1}^\Theta \Bigg[\sum_{i=1}^S \ee_{\ell n_\mu}  f_k(i) I_{\lbrace s(k) = i\rbrace} - F_{\min}(\check \theta) \Bigg] I_{\lbr \theta(\ell n_\mu) = \check \theta\rbr}\Bigg]\Bigg]\nonumber\\
& = \lim_{\mu\to 0}\ee h\left(\X^\mu(t_\iota),\theta^\mu(t_\iota):\iota\leq\kappa_0\right)\Bigg[ \sum^{ t+u}_{\ell:\ell \delta_\mu =t } \frac{\delta_\mu}{n_\mu} \nfm  \sum_{k=\ell n_\mu}^{\ell n_\mu+n_\mu-1} \nonumber\\
&\hspace{0.5cm}\Bigg[ \sum_{\check \theta = 1}^\Theta \Bigg[ \sum_{i=1}^S   \bbr_i\big(\fhat(\ell n_\mu)\big) \ee_{\ell n_\mu} f_k(i) - F_{\min}(\check \theta)\Bigg]   I_{\lbr \theta(\ell n_\mu) = \check \theta\rbr} \Bigg]\Bigg]\nonumber\\
&= \lim_{\mu\to 0}\ee h\left(\X^\mu(t_\iota),\theta^\mu(t_\iota):\iota\leq\kappa_0\right) \Bigg[ \sum^{ t+u}_{\ell:\ell \delta_\mu =t } \delta_\mu \nfm  \nonumber\\
& \hspace{0.7cm}\times\Bigg[ \frac{1}{n_\mu}\sum_{k=\ell n_\mu}^{\ell n_\mu+n_\mu-1} \sum_{\check \theta = 1}^\Theta I_{\lbr \theta^\mu(\mu\ell n_\mu) = \check \theta\rbr} \nonumber\\
&\hspace{1.2cm}\times\Bigg[ \sum_{i=1}^S   \bbr_i\big(\f^\mu(\mu\ell n_\mu)\big) \ee_{\ell n_\mu} f_k(i) - F_{\min}(\check \theta)\Bigg]\Bigg]\Bigg].
\end{align}
Here, we used $\ee_{\ell n_\mu} I_{\lbr s(k) = i\rbr} = \bbr_i\big(\f(\ell n_\mu)\big)$ as in Step 1. Recall that $\f(k) = \fhat(k) + \bf(\theta(k))$; see~(\ref{eq:f-hat}). By weak convergence of $\theta^\mu\cd$ to $\theta\cd$, the Skorohod representation, and using $\mu \ell n_\mu\to u$ and~(C1)--(C2), it can then be shown
\begin{align}
\label{eq:A-20}
& \sum_{\check \theta = 1}^\Theta  \frac{1}{n_\mu} \sum_{k=\ell n_\mu}^{\ell n_\mu + n_\mu -1} \bbr_i\big(\f^\mu(\mu\ell n_\mu)\big) \ee_{\ell n_\mu} f_k(i) I_{\lbr \theta^\mu(\mu\ell n_\mu) = \check \theta\rbr}  \nonumber\\
&\quad\to \bbr_i\big(\fhat(u) + \bf(\theta(u))\big) f(i,\theta(u))  \quad \textmd{in probability as } \mu \to 0.
\end{align}
Using a similar argument for the second term in~(\ref{eq:A1-19}), we conclude that, as $\mu\to 0$,
\begin{align}
\label{eq:A1-21}
&\ee h\left(\X^\mu(t_\iota),\theta^\mu(t_\iota):\iota\leq\kappa_0\right)\nonumber\\
 & \quad \hfill \times \Bigg[ \sum^{ t+u}_{\ell:\ell \delta_\mu =t } \delta_\mu \nfm \Bigg[\frac{1}{n_\mu} \sum_{k=\ell n_\mu}^{\ell n_\mu+n_\mu-1} \lb f_k(s(k)) - F_{\min}(\theta(\ell n_\mu))\rb\Bigg]\Bigg]\nonumber\\
& \to \ee h\left(\X^\mu(t_\iota),\theta^\mu(t_\iota):\iota\leq\kappa_0\right) \nonumber\\
&\hspace{0.5cm}\times\bigg[ \int_{t}^{t+u} \nfm(\X(v),\theta(v))\\
&\hspace{1cm}\times\lb \br\big(\fhat(v)+\bf(\theta(v))\big)\cdot \bf(\theta(v)) - F_{\min}(\theta(v))\rb dv\bigg].\nonumber
\end{align}
Finally, similar to (\ref{eq:CIR}),
\begin{align}
\label{eq:A1-22}
&
\ee h\big(\X^\mu(t_\iota), \theta^\mu(t_\iota): \iota\le \kappa_0\big) \nonumber\\
 &\quad\hfill \times\Bigg[\sum_{\ell:\ell \delta_\mu=t}^{t+u} \delta_\mu \nfm\Bigg[\frac{1}{n_\mu} \sum_{k = \ell n_\mu}^{\ell n_\mu+n_\mu-1} r(k)\Bigg]\Bigg]\nonumber\\
& \to \ee h\big(\X(t_\iota), \theta(t_\iota): \iota\le \kappa_0\big)\nonumber
\\
& \qquad\;\;\hfill \times \left[ \int^{t+u}_t
 \nfm(\X(v),\theta(v))\  r(v) dv\right]\;\;\textmd{as}\;\mu\to 0.
\end{align}

Combining the above two steps concludes the proof.

\section{Proof of Theorem~\ref{theorem:Stability}}
\label{sec:proof_2}

We first prove that each subsystem (the ODE~(\ref{eq:r-ode}) associated with each $\bar\theta\in\mathcal{Q}$ when $\theta(t) = \bar\theta$ is held fixed) is globally asymptotically stable $\mathbb{R}_{[0,\eta)}$ is its global attracting set. Define the Lyapunov function:
\begin{equation}
\label{eq:lyapunov_function}
V_{\bar\theta}\left(r\right)= r^2.\nonumber
\end{equation}
Taking the time derivative, and applying~(\ref{eq:r-ode}), we obtain
\begin{equation}
\label{eq:proof_B2}
\frac{d}{dt} V_{\bar\theta}\left(r\right) = 2r\cdot\lb \br\left(\bf(\bar \theta)\right)\cdot \bf(\bar\theta) - F_{\min}(\bar{\theta}) - r\rb\nonumber
\end{equation}
Since the objective function value at various states is bounded for each $\bar\theta\in\mathcal{Q}$,
\begin{equation}
\label{eq:proof_B3}
\frac{d}{dt} V_{\bar\theta}\left(r\right) \leq 2r\cdot\lb C(\gamma,\bar\theta) - r\rb\nonumber
\end{equation}
for some constant $C(\gamma,\bar\theta)$. Recall the smooth best-response sampling strategy $\br(\cdot)$ in Definition~\ref{def:smooth-best}. The parameter $\gamma$ simply determines the magnitude of perturbations applied to the objective function. It is then clear that $C(\gamma,\bar\theta)$ is monotonically increasing in $\gamma$.

In view of~(\ref{eq:proof_B3}), for each $\eta> 0$, $\hat{\gamma}$ can be chosen small enough such that, if $\gamma \leq \hat{\gamma}$ and $r \geq \eta$,
\begin{equation}
\label{eq:proof_B4}
\frac{d}{dt} V_{\bar\theta}\left(r\right) \leq - V_{\bar\theta}(r)\nonumber
\end{equation}
Therefore, each subsystem is globally asymptotically stable and, for $\gamma \leq \hat{\gamma}$,
\begin{equation}
\label{eq:proof_B5}
\lim_{t\to\infty} d\left(r(t),\RR_{[0,\eta)}\right) = 0.\nonumber
\end{equation}

Finally, stability of the regime-switching ODE~(\ref{eq:ODE-global}) is examined.
We can use the above
Lyapunov function to extend \cite[Corollary 12]{chatterjee2007switch} to prove global asymptotic stability w.p.1.

\vspace{0.1cm}
\begin{Theorem}[\cite{chatterjee2007switch}, Corollary 12]
\label{theorem:switched}
Consider the switching system~(\ref{eq:switched system}) in Definition~\ref{def:stability}, where $\theta(t)$ is the state of a continuous time Markov chain with generator $Q$. Define $\bar{q} := \max_{\theta\in\mathcal{Q}} |q_{\theta\theta}|$ and $\tilde{q} := \max_{\theta,\theta^{\prime}\in\mathcal{Q}} q_{\theta\theta^{\prime}}$. Suppose there exist continuously differentiable functions $V_{\theta}:\mathbb{R}^n\rightarrow\mathbb{R}^+$, $\theta\in\mathcal{Q}$, strictly increasing functions $a_1,a_2:\mathbb{R}^+\rightarrow\mathbb{R}^+$ with $a_1(0) = a_2(0) = 0$ and $a_1(t),a_2(t)\rightarrow\infty$ as $t\rightarrow\infty$, a real number $v>1$ such that the following hold:

\begin{enumerate}
    \item $a_1(d(Y,\H))\leq V_\theta(Y)\leq a_2(d(Y,\H)),\;\forall Y\in\RR^r,\theta\in\mathcal{Q}$,
    \item $\frac{\partial V_\theta}{\partial X} f_\theta(Y) \leq -\lambda V_\theta(Y)$, $\forall Y\in\RR^r,\forall \theta\in\mathcal{Q}$,
    \item $V_\theta(Y)\leq v V_{\theta^\prime}(Y)$, $\forall Y\in\RR^r, \theta,\theta^\prime\in\mathcal{Q}$,
    \item $(\lambda + \tilde{q})/\bar{q} > v$.
\end{enumerate}

\noindent
Then, the regime-switching system~(\ref{eq:switched system}) is globally asymptotically stable almost surely.
\end{Theorem}
\vspace{0.1cm}

The quadratic Lyapunov functions~(\ref{eq:lyapunov_function}) satisfies Hypothesis~2) in Theorem~\ref{theorem:switched}; see~(\ref{eq:proof_B4}). Further, since the Lyapunov functions are the same for all subsystems $\theta\in\mathcal{Q}$, existence of $v>1$ in Hypothesis~3) is automatically guaranteed. Hypothesis 4) simply ensures that the switching signal $\theta(t)$ is slow enough. Given that $\lambda = 1$ in hypothesis~2), it remains to ensure that the generator $Q$ of Markov chain $\theta(t)$ satisfies $1+\tilde{q} > \bar{q}$. This is satisfied since $|q_{\theta\theta^{\prime}}|\leq 1$ for all $\theta,\theta^{\prime}\in\mathcal{Q}$; see~(\ref{eq:Markov_cont}).

\bibliographystyle{IEEEtran}
\bibliography{IEEEabrv,DSO}

\newpage

\vfill

\end{document}